\documentclass[12pt]{article}
\usepackage{fullpage}
\usepackage{amsmath}
\usepackage{amssymb}
\usepackage{amscd}
\usepackage{epsfig}
\usepackage{rotating}

\begin{document}
\input{amssym.def}

\title{Cohomology of line bundles on Schubert varieties in the Kac-Moody setting}\author{S.~Senthamarai Kannan}

\date{}
\def\bp{{\mathbb P}}
\def\bc{{\mathbb{C}}}
\newcommand{\lr}{\longrightarrow}
\newcommand{\ctext}[1]{\makebox(0,0){#1}}

\newcommand{\globalg}[1]{\mbox{\rm $H^0(G/B,#1)$}}
\newcommand{\globalh}[1]{\mbox{\rm $H^0(H/B \cap H,#1)$}}
\newcommand{\cohom}[3]{\mbox{\rm $H^{#1}(#2,#3)$}}
\newcommand{\lprod}[2]{\mbox{\rm $\langle #1~,{{#2}}\rangle$}}
\newcommand{\lb}[1]{\mbox{\rm $L_{#1}$}}
\newcommand{\high}[1]{\mbox{\rm $(x_1 + x_{2n +1})^{#1}$}}
\newcommand{\sln}{\mbox{\rm SL$(2n+1,K)$}}
\newcommand{\son}{\mbox{\rm SO$(2n+1,K)$}}
\newcommand{\gsec}[1]{\mbox{\rm $s_{#1}$}}
\newcommand{\gsecb}[1]{\mbox{$\overline{s_{#1}}$}}
\newcommand{\chevzform}[1]{\mbox {$\mathcal{U}_{\mathbb{Z}}({#1})$}}
\newcommand{\chevqform}[1]{\mbox {$\mathcal{U}_{\mathbb{Q}}({#1})$}}
\newcommand{\roots}{R}
\newcommand{\proots}{R^{+}}
\newcommand{\nroots}{R^{-}}
\newcommand{\lie}[1]{\mbox {\cal{#1}}}
\newcommand{\uzthspan}{\left<\chevzform{\lie{h}}v^{+}\right>}
\newcommand{\univenv}[1]{\mathcal{U}(\cal{#1})}
\newcommand{\uzchspan}{\left<\chevzform{\lie{h}}v^{+}\right}
\newcommand{\Pow}[1]{{\it power}_{#1}\hbox{\rm-GA}}
\newcommand{\A}[1]{{\it Aut}(#1)}
\newcommand{\id}{{\rm id}}
\newcommand{\path}{\hbox{Path}}
\newcommand{\cycle}{\hbox{Cycle}}
\newcommand{\Pm}{\leq_m^p}
\newcommand{\PT}{\leq_T^p}
\newcommand{\ET}{\equiv_T^p}
\newcommand{\pmequiv}{\equiv_m^p}
\newcommand{\lequivr}{\leftrightarrow}
\newcommand{\limp}{\rightarrow}
\newcommand{\perm}[2]{ \sigma_{#1,#2}}
\newcommand{\zmod}[1]{\mu_{#1}}
\newcommand{\jn}[1]{~{\bf #1}}
\newcommand{\qed}{\hfill$\Box$}
\newcommand{\irredH}[1]{L(#1)_K}

\newcommand{\beqn}{\begin{eqnarray*}}
\newcommand{\eeqn}{\end{eqnarray*}}
\newcommand{\pair}[2]{\mbox {$(#1,#2)$}}
\newcommand{\opone}[1]{\mbox{$Y_{\alpha_{#1}}$}}
\newcommand{\optwo}[2]{\mbox{$Y_{\alpha_{#1}}Y_{\alpha_{#2}}$}}
\newcommand{\opthree}[3]{\mbox{$Y_{\alpha_{#1}}Y_{\alpha_{#2}}Y_{\alpha_{#3}}$}}
\newcommand{\opfour}[4]{\mbox{$Y_{\alpha_{#1}}Y_{\alpha_{#2}}Y_{\alpha_{#3}}Y_{\alpha_{#4}}$}}

\newtheorem{guess}{Theorem}[section]
\newcommand{\bth}{\begin{guess}$\!\!\!${\bf .}~~\rm}
\newcommand{\eeth}{\end{guess}}

\newtheorem{mtheo}[guess]{\bf Conjecture}
\newcommand{\bbth}{\begin{mtheo}$\!\!\!${\bf .}~~\rm}
\newcommand{\eeeth}{\end{mtheo}}

\newtheorem{propo}[guess]{Proposition}
\newcommand{\bprop}{\begin{propo}$\!\!\!${\bf .}~~\rm}
\newcommand{\eprop}{\end{propo}}

\newtheorem{equo}[guess]{Equation}
\newcommand{\bq}{\begin{equation}}
\newcommand{\eq}{\end{equation}}

\newtheorem{claim}[guess]{\it Claim}
\newcommand{\bcl}{\begin{claim}$\!\!\!${\bf .}~~\rm}
\newcommand{\ecl}{\end{claim}}

\newtheorem{rema}[guess]{\it Remark}
\newcommand{\brem}{\begin{rema}$\!\!\!${\it .}~~\rm}
\newcommand{\erem}{\end{rema}}
\newtheorem{coro}[guess]{Corollary}
\newcommand{\bcor}{\begin{coro}$\!\!\!${\bf .}~~\rm}
\newcommand{\ecor}{\end{coro}}
\newtheorem{lema}[guess]{Lemma}
\newcommand{\blem}{\begin{lema}$\!\!\!${\bf .}~~\rm}
\newcommand{\elem}{\end{lema}}
\newtheorem{dem}[guess]{Demazure Sequence}
\newcommand{\bdem}{\begin{dem}$\!\!\!${\bf .}~~\rm}
\newcommand{\edem}{\end{dem}}

\newtheorem{exam}{Example}
\newcommand{\beg}{\begin{exam}$\!\!\!${\bf .}~~\rm}
\newcommand{\eeg}{\end{exam}}
\newcommand{\er}{\hfill {\Large $\bullet$}\linebreak}
\newtheorem{defe}[guess]{Definition}
\newcommand{\bdefe}{\begin{defe}$\!\!\!${\bf .}~~ \rm}
\newcommand{\edefe}{\end{defe}}
\newcommand{\spec}{{\rm Spec}\,}

\newcommand{\mybar}[1]{\mbox{\overline{#1}}}
\newtheorem{theorem}{Theorem} 
\newtheorem{lemma}[theorem]{Lemma}
\newtheorem{proposition}[theorem]{Proposition}
\newtheorem{corollary}[theorem]{Corollary}
\newtheorem{definition}[theorem]{Definition}
\newtheorem{observation}[theorem]{Observation}
\newtheorem{question}[theorem]{Question}

\newenvironment{proof}{\noindent{\it Proof.  }}{\hspace*{\fill}
\rule{2mm}{2mm} \vspace{5mm}}

\newcommand{\claimproof}[2]%
{\noindent{\it Proof of Claim#1.}
#2\hspace*{\fill}$\Box$~~~~~\vspace{5mm} } 

\maketitle

\begin{abstract}In this paper, we describe the indices of the top and the least nonvanihing  cohomologies $H^{i}(X(w), L_{\lambda})$ of line budles on Schubert varieties $X(w)$ given by nondominant weights in the Kac-Moody setting. We also prove some surjective Theorem for maps between some cohomology modules.
\end{abstract}  

\section{Introduction}

The aim of this paper is to study the cohomology modules of line bundles on Schubert varieties given by non-dominant weights in the {\it Kac-Moody} setting.  

The following notations will be maintained throughout this paper.
The base field is $\bc$, the field of complex numbers.
Let $\mathcal{G}$ be a Kac Moody Lie algebra associated to a generalised symmetrizable Cartan matrix [cf [7]]. 

Let $\mathcal{G}$ be a Kac-Moody group associated to 
$\mathcal{G}$ [ cf pp. 183[9]]. Let $\mathcal{H}$ be a Cartan 
subalgebra of $\mathcal{G}$.
Let $T$ be a maximal torus of $\mathcal{G}$ with $\mathcal{H}$ as its lie algebra [cf pp.178, 8] and let $X(T)$ denote the set of characters of $T$, $W$ denote the Weyl group of $\mathcal{G}$ with respect to $T$. Let $R$ denote the set of roots of $G$ with respect to $T$, $\mathcal{B}=T\mathcal{U}$ be a Borel subgroup of $\mathcal{G}$ [ cf pp.175, pp.183 [8]]. Let $S=\{\alpha_1,\ldots,\alpha_l\}$ denote the set of simple roots in $\proots$. For $\beta \in R^+$ we also use the notation $\beta > 0$. The element of the Weyl group (i.e. the simple reflection) corresponding to $\alpha_i$ is denoted by $s_{\alpha_i}$. For a simple root $\gamma$, we denote the smallest subgroup of $B$ containing $T$ and the root subgroup $G_{a, -\gamma}$ by $B_{\gamma}$. A non degenerate $W$-invariant bilinear form on $X(T)$ induced by the one on $\mathcal{H}$ [cf pp.17[7]] is denoted by $(~,~)$. 
When $\alpha$ is a real root, then we use the following notation $\left<~,~ \right>$ to denote $\lprod{\nu}{\alpha} = \frac{2(\nu, \alpha)}{(\alpha,\alpha)}$.
Note that $(\alpha, \alpha)$ is positive when $\alpha$ is a real root [cf [7]].
 We denote by $\Lambda^+$ the set of dominant weights i.e. the set of weights $\lambda \in\Lambda$, such that $\lprod{\lambda}{\alpha} \geq 0$ for all $\alpha
\in R^+$. We denote by $\Lambda^{++}$ the set of regular dominant weights. We fix an element $\rho\in \Lambda^+$ such that $\lprod{\rho}{\alpha_{i}}=1$ for every $i=1, 2, \cdots l$. 

For $w \in W$ let $l(w)$ denote the length of $w$. For $w \in W$, let $X(w)$ denote the Schubert variety in $\mathcal{G}/\mathcal{B}$ corresponding to $w$. Let $\leq$ denote the Bruhat order on $W$.  

When $\mathcal{G}$ is finite dimensional, a systematic study of the cohomology modules of line bundles on Schubert varieties given by non-dominant weights was done in \cite{bwbk}. In this paper, we undertake a systematic study of the cohomology modules of line bundles on Schubert varieties given by non-dominant weights in the Kac-Moody setting. Besides, several results in this paper are new, even in the finite dimensional case. 

Broadly, the paper is in the same spirit as in \cite{bwbk}. As in \cite{bwbk}, the strategy is a delicate use of the Bott-Samelson inductive machinary, available to us in the Kac-Moody setting from \cite{kum2}, \cite{Math} and \cite{Slod}. However, we would like to point out that there is a technical difficulty in carrying out the proofs as in \cite{bwbk}. The proofs in \cite{bwbk} are based on descending induction on the dimension of the Schubert variety using the Borel-Weil-Bott Theorem for the flag variety $G/B$. This argument would not be available in the Kac-Moody setting, since the flag variety is infinite dimensional even though, the Borel-Weil-Bott Theorem holds when $G$ is infinite dimensional (see below). Furthermore, a number of the proofs in \cite{bwbk} use descending induction on the length of the element in the Weyl group moving the non-dominant weight into the dominant chamber. And so the proofs in this paper are more subtle. 

In the case when $\lambda$ is dominant, this problem has been well
studied even in the Kac-Moody setting. When $\lambda$ is non-dominant 
and when $X(w)\simeq G/B$ and $G$ is not finite dimensional, the 
equivalent of the Borel-Weil-Bott Theorem is known to be true  and was proved 
independently by S. Kumar (cf \cite{kum} ) and O. Mathieu (cf \cite{Math}). To the best of our knowledge, the cohomology of line bundles associated to non-dominant weights on Schubert varieties in the Kac-Moody setting has never been addressed. And, this is the focus of the current paper.

For a uniform theory to hold, we need to make certain {\it genericity} assumptions on the weights. In all the situations, this has been made very specific to the problem at hand. The case when the weights are somewhat {\it special}, situated essentially near the walls of some Weyl chambers, the behaviour can be erratic and seems to involve very complicated combinatorics. 

We consistently use the following terminologies. Consider the Tits Cone $X=\bigcup_{w\in W} w(\Lambda^{++})-\rho$. Let $\phi\in W$. A weight $\lambda\in X$ such that $\phi\cdot\lambda$ is dominant is said to be generic if for all simple roots $\alpha$, one has $|\lprod{\lambda}{\alpha}| \gg 0$. For a generic weight $\lambda$, it is clear that the element $\phi$ is unique. We then say that $\lambda$ is a generic weight in the $\phi$-chamber. Let ${\mathcal L}_{\lambda}$ denote the line bundle on $X(w)$ corresponding to the 1-dimensional representation of $B$ given by the character $\lambda$. In such a situation, for the cohomology module $\cohom{i}{X(w)}{{\cal L}_{\lambda}}$ we simply write $\cohom{i}{w}{\lambda}$. Wherever needed, we have indicated the precise genericity conditions.

The layout of this paper is as follows:

In section 2, we setup our notation and recall some basic Theorems needed. The standard results on the combinatorics of the Weyl group that are used in this paper can be found in \cite {kac}. In section 3, we state  one combinatorial Lemma and some useful corollaries involving the elements of the Weyl group and weights.  

In section 4, we prove the following Theorem. This Theorem is known when $G$ is finite dimensional [cf \cite{bwbk}, \cite{Dab}, \cite{Polo}]: 

{\bf Theorem} Let $X(w)$ be a Schubert variety and $\lambda$ a
generic weight in the $\phi$-chamber. Let $R^+(w) = \{\alpha \in R^+ \mid
w(\alpha) \in R^-\}$.
\begin{enumerate}
\item $\cohom{0}{w}{\lambda} \neq 0$ if and only if $R^+(w) \cap
R^+(\phi) = \emptyset$.
\item $\cohom{l(w)}{w}{\lambda} \neq 0$ if and only if $R^+(w) \subseteq R^+(\phi)$.
\item The restriction map $H^{l(\phi)}(G/B,\lambda)\longrightarrow H^{l(\phi)}(\phi,\lambda)$ is surjective.
\end{enumerate}

In section 5, we prove the following Theorem. Note that the statement $(2)$ of this Theorem is new, even in the finite dimensional case. The statements $(1)$ and $(3)$ of this Theorem are known to be true in the finite dimensional case [cf \cite{bwbk}, \cite{Jo}]:

{\bf Theorem} Let $X(w)$ be a Schubert variety and $\lambda$ be a
generic weight in the $\phi$-chamber. 
\begin{enumerate}
\item $\cohom{j}{w}{\lambda} = 0$ for $j > \min(l(w),l(\phi))$.
\item If $\phi\leq w$ , then $\cohom{l(\phi)}{w}{\lambda}\neq 0$ and the restriction map $H^{l(\phi)}(G/B, \lambda)\longrightarrow H^{l(\phi)}(w,\lambda)$ is surjective. 
\item (Cohomological characterization of the Bruhat order)

$\phi \leq w$ if and only if $\cohom{l(\phi)}{w}{\lambda} \neq 0$.
\end{enumerate}
To describe the results in section 6, we first introduce the following notation
[cf section 6]. Let $W^{+}(w,\phi):=\{\tau\leq w : R^{+}(\tau)\subset R^{+}(\phi)\}$ and let $W^{-}(w,\phi):=\{\tau\leq w: R^{+}(\tau)\bigcap R^{+}(\phi) ~ is ~ empty\}$. We prove that both the sets have unique maximal elements with respect to the Bruhat order, which we denote by $\tau^{+}(w,\phi)$ and $\tau^{-}(w,\phi)$ respectively. We call $l^{+}(w,\phi)=l(\tau^{+}(w,\phi))$ and $l^{-}(w,\phi)=l(\tau^{-}(w,\phi))$. Let $M:=\max\{\langle \beta, \gamma\rangle: \beta\in \phi(S), \gamma\in S\}$.      

In section 6, we prove the following Theorem:

{\bf Theorem} Let $X(w)$ be a Schubert variety and let $\lambda$ a generic weight in the $\phi$-chamber. Let $D$ denote the boundary divisor of $X(w)$.
\begin{enumerate}
\item If $\lambda$ is a weight such that $\phi\cdot\lambda$ is dominant,
then, the restriction map $H^{l^{+}(w,\phi)}(w,\lambda)
\longrightarrow H^{l^{+}(w,\phi)}(\tau^{+}(w,\phi),\lambda)$ is surjective 
and in particular $H^{l^{+}(w,\phi)}(w,\lambda)$ is a non zero $B$- module 
with $\tau^{+}(w,\phi)\cdot\lambda$ as a weight of this $B$- module. 
\item When $G$ is finite dimensional, and if $\langle \phi\cdot\lambda , \gamma \rangle >  l(w _{0}\phi)M$ (here $w_{0}$ denote the longest element of the Weyl group of $G$), for all simple roots $\gamma$, then, we have $H^{i}(w, \lambda)=0$ for $i > l^{+}(w, \phi)$.
\item If $\lambda$ is a weight such that $\langle \phi\cdot\lambda, \gamma \rangle ~ > ~ l(\phi)M$ for all simple roots $\gamma$, then, we have $H^{i}(w,\lambda)=0$ for all $i ~ < ~ l(w)-l^{-}(w,\phi)$.
\item If $\lambda$ is a weight such that $\langle \phi(\lambda), \gamma \rangle \geq -1$ for all simple roots $\gamma$, then, $H^{l(w)-l^{-}(w, \phi)}(w,\lambda- D)\neq 0$.
\end{enumerate}

\section{Preliminaries}
\label{prelim}

Given a $w \in W$ the closure in $G/B$ of the $B$ orbit of the coset
$wB$ is the Schubert variety corresponding to $w$, and is denoted by
$X(w)$. We recall some basic facts and results about Schubert
varieties. A good reference for all this is Jantzen's book [cf
~\cite[II, Chapter 14 ]{J}] in the finite dimensional case,  [cf \cite[Sections 7.1 and  8.1]{kum2}], [cf \cite{Math}] and [cf \cite{Slod}] in the infinite dimensional 
case.

For a simple root $\alpha$, we denote by $P_{\alpha}$ the minimal parabolic subgroup of $G$ containing $B$ and $s_{\alpha}$. If $\alpha=\alpha_{i}$, we also denote the minimal parabolic subgroup $P_{\alpha_{i}}$ by $P_{i}$.  
Let $w = s_{\alpha_{i_{1}}}s_{\alpha_{i_{2}}}\ldots s_{\alpha_{i_{n}}}$ be a reduced expression for $w \in W$. Define 
\[Z(w) = \frac {P_{i_{1}} \times P_{i_{2}} \times \ldots \times P_{i_{n}}}{B \times \ldots
\times B}
\]
where the action of $B \times \ldots \times B$ on $P_{i_{1}} \times P_{i_{2}}
\times \ldots \times P_{i_{n}}$ is given by $(p_1, \ldots , p_n)(b_1, \ldots
, b_n) = (p_1 \cdot b_1, b_1^{-1} \cdot p_2 \cdot b_2, \ldots
,b^{-1}_{n-1} \cdot p_n \cdot b_n)$, $ p_i \in P_i$, $b_i \in B$. Note
that $Z(w)$ depends on the reduced expression chosen for $w$. It is
well known that $Z(w)$ is a smooth projective $B$- variety and there exists a birational surjective morphism
\[
\phi_w : Z(w) \lr X(w)
\]

Let $f_n : Z(w) \lr Z(ws_{\alpha_{i_{n}}})$ denote the map induced by the
projection $P_{i_{1}} \times P_{i_{2}} \times \ldots \times P_{i_{n}} \lr 
P_{i_{1}} \times P_{i_{2}}\times \ldots \times P_{i_{n-1}}$. We observe that $f_n$ is a $ P_n/B\simeq {\bf P}^{1}$-fibration. We also denote ${\bf P}^{1}=P_{\alpha}/B$ by ${\bf P}_{\alpha}^{1}$.  

Let $V$ be a $B$-module. Then, for each $w \in W$, we obtain by the
standard method of associated construction an induced bundle ${\cal
L}_w(V)$ on $X(w)$ and then on $Z(w)$ via the map $\phi_{w}$. Then, for $i \geq 0$ we have the following isomorphisms of $B$-linearized sheaves
\[
R^{i}{f_{n}}_{*}{\mathcal L}_w(V) = {\mathcal
L}_{ws_{\alpha_{i_{n}}}}(\cohom{i}{P_{i_{n}}/B}{{\mathcal L}_w(V)})
\]
This together with easy applications of Leray spectral sequences is
the constantly used tool in what follows. We term this the {\em
descending 1-step construction}. 

We also have the {\em ascending 1-step construction} which is used
extensively in what follows in conjunction with the descending construction. We recall this for the convenience of the reader.

Let the notation be as above and write $\tau = s_{\alpha}w$, with 
$l(\tau) = l(w) +1$, for some simple root $\alpha$. Then, we have an induced morphism
\[
g_{\alpha}: Z(\tau) \lr P_{\alpha}/B \simeq {\bf P}^1
\]
with fibres given by $Z(w)$. Again, by an application of the Leray
spectral sequences together with the fact that the base is a ${\bf
P}^1$, we obtain for every $B$-module $V$ the following exact sequence
of $P_{\alpha}$-modules
\[
0 \lr \cohom{1}{P_{\alpha}/B}{R^{i-1}{g_{\alpha}}_{*}{\mathcal L}_w(V)} \lr
\cohom{i}{Z(\tau)}{{\mathcal L}_{\tau}(V)} \lr
\cohom{0}{P_{\alpha}/B}{R^{i}{g_{\alpha}}_{*}{\mathcal L}_w(V)} \lr 0
\]

We also recall the following well-known isomorphisms [cf II Chapter 14 in \cite{J} in fin. diml., and [ cf 7.1 and 8.1 in \cite{kum2} in inf. diml. case]:
\begin{itemize}

\item ${\phi_w}_*{\cal O}_{Z(w)} = {\cal O}_{X(w)}$.

\item $R^{q}{\phi_w}_*{\cal O}_{Z(w)} = 0$ for $q > 0$.

\end{itemize}

This together with \cite[II. 14.6]{J} implies that we may use the
Bott-Samelson schemes $Z(w)$ for the computation and study of all the
cohomology modules $\cohom{i}{X(w)}{{\mathcal L}_w(V)}$. Henceforth in this paper, we will use the Bott-Samelson schemes and their cohomologies in all the computations.

\noindent
{\bf Notation}: 

{\em Here and in what follows we have by an abuse of notation replaced
the induced sheaf ${\mathcal L}(\cohom{i}{P_{\alpha}/B}{{\mathcal
L}_{\lambda}})$ by the cohomology module $\cohom{i}{\alpha}{\lambda}$
and the cohomology modules $\cohom{i}{Z(w)}{{\mathcal L}_{V}}$ by the
$\cohom{i}{w}{V}$}.

\subsubsection {Some constructions from Demazure's paper}
We recall briefly two exact sequences that Demazure used in his short
proof of the Borel-Weil-Bott theorem~\cite{D}. We use the same 
notation as in Demazure. In the rest of the paper these sequences are 
referred to as Demazure exact sequences. We first recall the ``dot'' action of the Weyl group $W$ on $X(T)$. Let $w\in W$ and $\lambda\in X(T)$. For $\lambda\in X(T)$, we define $w\cdot\lambda:=w(\lambda+\rho)-\rho$.  

Let $\alpha$ be a simple root and let $\lambda \in X(T)$ be a weight
such that $\lprod{\lambda}{\alpha} \geq 2$. Set $V_{\lambda, \alpha}:=H^{0}(P_{\alpha}/B, {\mathcal L}_{\lambda})$. It is easy to see that $V_{\lambda,\alpha}$ is an irreducible $P_{\alpha}$- module and it breaks up as a $T$-module into a direct sum of weight spaces with weights $\lambda,\lambda-\alpha,\ldots, s_{\alpha}(\lambda)$. Let $L_{\lambda}$ denote the 1-dimensional representation given by the character $\lambda$ of $B$ then as $B$-modules we have the following exact sequences:

\bdem
\label{Dem 1}
$$
\begin{array}{l}
\mbox{$0 \lr K \lr V_{\lambda,\alpha} \lr L_\lambda \lr 0$}\\
\mbox{$0 \lr L_{s_{\alpha}(\lambda)} \lr K \lr V_{\lambda-\alpha,\alpha}
\lr 0$}\\
\end{array}
$$
\edem

Dually, if $\langle \lambda, \alpha\rangle \leq -2$, then, we have:

\label{Dem 2}
$$
\begin{array}{l}
\mbox{$0 \lr L_{\lambda} \lr V_{s_{\alpha}(\lambda),\alpha} \lr Q \lr 0$}\\
\mbox{$0 \lr V_{s_{\alpha}\cdot\lambda, \alpha} \lr Q \lr L_{s_{\alpha}(\lambda)} \lr 0$}\\
\end{array}$$

A consequence of the above exact sequences is the following lemma. The
proof is exactly as in Demazure but we give it here for completeness.
\blem \label{Spec-1}
\begin{enumerate}
\item Let $\tau = ws_{\alpha}$, $l(\tau) = l(w)+1$. If
$\lprod{\lambda}{\alpha} \geq 0$ then for all $j$,
$\cohom{j}{\tau}{\lambda} = H^{j}(w, V_{\lambda, \alpha})$.
\item Let $\tau = ws_{\alpha}$, $l(\tau) = l(w)+1$. If
$\lprod{\lambda}{\alpha} \geq 0$ then $\cohom{i}{\tau}{\lambda} =
\cohom{i+1}{\tau} {s_{\alpha}\cdot \lambda}$ and if
$\lprod{\lambda}{\alpha} \leq -2$ then $\cohom{i}{\tau}{\lambda} =
\cohom{i-1}{\tau}{s_{\alpha}\cdot \lambda}$
\item If $\langle \lambda, \alpha \rangle =-1$, then, $H^{i}(\tau,\lambda)=(0)$
for all $i$.
\end{enumerate}
\elem
\begin{proof}
\begin{enumerate}
\item Let $f: Z(\tau) \lr Z(w)$ be the ${\bf
P}^{1}_{\alpha}$-fibration. Then we have
\[
f_{*}{\mathcal L}_{\lambda} = {\mathcal
L}(\cohom{0}{P_{\alpha}/B}{{\mathcal L}_{\lambda}})
\]
Therefore, by an application of Leray spectral sequences (as in \cite[II 14.6]{J} )it follows
that $\cohom{j}{\tau}{\lambda} = H^{j}(w, V_{\lambda, \alpha})$,
since $\cohom{0}{P_{\alpha}/B}{{\mathcal L}_{\lambda}} =
V_{\lambda, \alpha}$ and since  
\[
R^{1}f_{*}{\mathcal L}_{\lambda} = {\mathcal
L}(\cohom{1}{P_{\alpha}/B}{{\mathcal L}_{\lambda}}) = 0
\]
when $\lprod{\lambda}{\alpha} \geq 0$.

\item We prove it in the case $\lprod{\lambda}{\alpha} \geq 0$. Then
$H^1(\alpha, \lambda) = 0$.  Hence we get $\cohom{i}{\tau}{\lambda} =
H^i(w,H^0(\alpha,\lambda))$. But $H^0(\alpha, \lambda) = H^1(\alpha,
s_{\alpha} \cdot \lambda) $. Hence we get $\cohom{i}{\tau}{\lambda} =
H^i(w,H^1(\alpha,s_{\alpha}\cdot\lambda))$.  But this last module is
exactly $H^{i+1}(\tau, s_{\alpha}\cdot\lambda)$.  The other case is
similar.
\end{enumerate}
\end{proof}

\section{Combinatorial Lemmas:}

In this section, we prove a combinatorial Lemma and state some Corollaries of this Lemma.

\blem Let $w\in W$. Let $\lambda$ be any weight. Then, for any $i\in\{0,1, 2,\cdots l(w)\}$, every weight of $H^{i}(w,\lambda)$ is in the convex hull of $\{\tau\cdot\lambda:\tau\leq w\}$. \elem

\begin{proof} Proof is by induction on $l(w)$, the base case $w=id$   
being trivial. Let $l(w)$ be positive, and let $w=s_{\gamma}\tau$ with 
$\gamma$ a simple root and $l(w)=l(\tau)+1$. Now, let $i\in\{0, 1, 2, \cdots l(w)\}$, and let $\mu$ be a weight of $H^{i}(w,\lambda)$. Since, we have 
$(0)\longrightarrow H^{1}(s_{\gamma}, H^{i-1}(\tau,\lambda))\longrightarrow 
H^{i}(w, \lambda)\longrightarrow H^{0}(s_{\gamma}, H^{i}(\tau,\lambda))
\longrightarrow (0)$, $\mu$ must be of the form 
$\mu=a\nu+(1-a)s_{\gamma}\cdot\nu$ for some weight $\nu$ of $H^{i-1}(\tau,\lambda)$ or of $H^{i}(\tau, \lambda)$ and for some $0\leq a\leq 1$. By induction, 
$\nu$ is of the form $\nu= \sum_{\tau^{\prime}\leq \tau} c_{\tau^{\prime}} \tau^{\prime}\cdot \lambda$, with $0\leq c_{\tau^{\prime}}\leq 1$ for each 
$\tau^{\prime}$ and the sum $\sum_{\tau^{\prime}}c_{\tau^{\prime}}=1$.
Observe that $s_{\gamma}\cdot \sum_{\tau^{\prime}\leq \tau}
c_{\tau^{\prime}}\tau^{\prime}\cdot \lambda=\sum_{\tau^{\prime}\leq \tau}
c_{\tau^{\prime}}s_{\gamma}\tau^{\prime}\cdot \lambda$. Since 
$w^{-1}(\gamma)$ is a negative root, $s_{\gamma}\tau^{\prime}\leq w$ for any $\tau^{\prime}\leq \tau$. Hence the Lemma follows from the expression of $\mu$.
\end{proof}

We have:
\bcor  Let $\lambda$ be a weight such that $\phi\cdot\lambda$ is dominant.
Let $w, ~\tau\in W$ be arbitrary. If $\tau\cdot\lambda$ is a weight of 
$H^{i}(w, \lambda)$ for some $i$, then, $\tau\leq w$. \ecor
\begin{proof} Since $\lambda+\rho$ is non-singular, there is a $\phi\in W$ such that $\phi\cdot\lambda$ is dominant.  Now, if $\tau\cdot\lambda$ is a weight of $H^{i}(w, \lambda)$ for some $i$, then, by Lemma 3.1, there are real numbers $0\leq c_{\tau^{\prime}}\leq 1$, $\tau^{\prime}\leq w$ with $\sum_{\tau^{\prime}\leq w}c_{\tau^{\prime}}=1$ such that $\tau\cdot \lambda = \sum_{\tau^{\prime}\leq w}c_{\tau^{\prime}}\tau^{\prime}\cdot \lambda$. 
Hence, we must have $\phi\cdot\lambda=\sum_{\tau^{\prime}\leq w}
c_{\tau^{\prime}}(\phi \tau^{-1}\tau^{\prime}\phi^{-1})\cdot \phi\cdot\lambda$.

Now, since $\phi\cdot \lambda$ is dominant, the above equality holds only if 
$\tau=\tau^{\prime}$ for any $\tau^{\prime}\leq w$ such that 
$c_{\tau^{\prime}}\neq 0$. This forces that $\tau\leq w$.     
\end{proof}

\bcor  Let $\lambda$ be a weight such that $\phi\cdot\lambda$ is a dominant weight for some $\phi\in W$. Let $w\in W$ be arbitrary. Let $\tau\in W$, let $\gamma$  be a simple root such that $\phi\tau^{-1}(\gamma)$ is a positive root. 
If $\tau\cdot\lambda+m\gamma$ is a weight of $H^{i}(w, \lambda)$ for some $i$ and for some non-negative integer $m$, then, $\tau\leq w$. \ecor
\begin{proof} If $\tau\cdot\lambda+m\gamma$ is a weight of $H^{i}(w, \lambda)$ for some $i$, then, it follows from Lemma 3.1 that there are real numbers $0\leq c_{\tau^{\prime}}\leq 1$, $\tau^{\prime}\leq w$ with $\sum_{\tau^{\prime}\leq w}
c_{\tau^{\prime}}=1$ such that $\phi\cdot\lambda + m 
\phi\tau^{-1}(\gamma)=\sum_{\tau^{\prime}\leq w}
c_{\tau^{\prime}}(\phi \tau^{-1}\tau^{\prime}\phi^{-1})\cdot\phi\cdot\lambda$.
Now, since $\phi\tau^{-1}(\gamma)$ is a positive root, $m$ must be zero and 
$\tau\leq w$.  
\end{proof}

\section{Non vanishing of the zeroth and the topmost cohomology- neccessary
and sufficient conditions:}
In this section, we give a criterion for the non-vanishing of $H^{0}(w,\lambda)$ and $H^{l(w)}(w,\lambda)$. The conditions are based on the combinatorics between $w$ and $\phi$ if $\lambda$ is a generic weight in the $\phi$- chamber. 

Let $w$, $\phi$ be two elements of the Weyl group. Let $\lambda$ be a  
weight such that $\phi(\lambda)$ is regular dominant. 
Then, we have the following:

\bprop $H^{0}(w,\lambda)\neq 0$ if and only if $R^{+}(w)\bigcap R^{+}(\phi)$ 
is empty. Further, if $H^{0}(w, \lambda)\neq 0$, then the dimension
of the subspace of $U$-invariant vectors in $H^{0}(w,\lambda)$ is one and 
it is spanned by a weight vector of weight $w(\lambda)$. \eprop 
\begin{proof} The proof given in Theorem 3.3(i) of \cite{bwbk} for the `if' part holds also in the Kac-Moody setting. So, we need to prove only the converse. We prove this by induction on $l(\phi)$.

If $l(\phi)=0$, there is nothing to prove. 
Let $\phi\in W$ be an element whose length is positive.
Then, there is a simple root $\alpha$ such that 
$\phi(\alpha)$ is a negative root. If $w(\alpha)$ is a negative root,
then, $H^{0}(w,\lambda)=H^{-1}(w,s_{\alpha}\cdot \lambda)=0$.
Otherwise, consider the following exact sequence of $B$- modules:
$(0)\longrightarrow \lambda \longrightarrow V_{s_{\alpha}(\lambda), \alpha}
\longrightarrow Q\longrightarrow (0)$.

From this, applying $H^{0}$, we get 
$(0)\longrightarrow H^{0}(w, \lambda)\neq (0)\longrightarrow 
H^{0}(w s_{\alpha}, s_{\alpha}(\lambda))$.
(Here we note that  $\langle s_{\alpha}(\lambda), \alpha \rangle  > 0$ since 
$\phi(\alpha)$ is a negative root.)       

Now, since $l(\phi s_{\alpha})=l(\phi)-1$, and  from the above observation
since $H^{0}(w, \lambda)\neq 0$, by induction, we get $R^{+}(w s_{\alpha}) 
\bigcap R^{+}(\phi s_{\alpha})$ is empty. Therefore, 
$R^{+}(w)\bigcap R^{+}(\phi)$ is empty. (For if $\beta\in R^{+}(w)\bigcap R^{+}(\phi)$, then, $\beta\neq \alpha$, which forces that $s_{\alpha}(\beta)\in R^{+}(ws_{\alpha})\bigcap R^{+}(\phi s_{\alpha})$, a contradiction to the above observation).  

The weight computation also can be seen along the same lines of the 
above argument. \end{proof}

Let $w$, $\phi\in W$ be two elements of the Weyl group.
Let $\lambda$ be a weight such that $\lambda + \rho$ is non singular and 
$\phi \cdot \lambda$ is dominant. Then, we have
\bprop  $H^{l(w)}(w, \lambda)\neq 0$ if and only if 
$R^{+}(w)\subset R^{+}(\phi)$. When $H^{l(w)}(w,\lambda)\neq 0$, it is a 
cyclic $B$- module generated by a weight vector of weight $w\cdot \lambda$.
\eprop
\begin{proof} The proof of the `if' part for the finite dimensional case given in Theorem 3.3(ii) of \cite{bwbk} holds also in the Kac-Moody setting. So, we only need to prove that if $H^{l(w)}(w,\lambda)\neq 0$, then, $R^{+}(w)\subset R^{+}(\phi)$. 

We prove this in two steps.

{\bf Step(1):}
We first prove that if $R^{+}(w)\subset R^{+}(\phi)$, then, 
$H^{l(w)}(w,\lambda)$ is a cyclic $B$- module generated by a weight vector 
of weight $w\cdot \lambda$. 

By Serre-duality, we have an ismorphism $H^{l(w)}(w, \lambda)^{*}\sim H^{0}(w, -\lambda-\rho-D)\bigotimes L_{\Psi_{w}}$ of $B$- modules for some character $\Psi_{w}$ of $B$. Here, $D$ denote the sheaf associated to the boundary divisor of $X(w)$, and $L_{\Psi_{w}}$ is the one dimensional $B$- module given by the character $\Psi_{w}$ of $B$. On the otherhand, we have an inclusion $H^{0}(w, -\lambda-\rho-D)\bigotimes L_{\chi_{w}}\hookrightarrow H^{0}(w,-\lambda-\rho)$ of $B$- modules for some character $\chi_{w}$ of $B$. Therefore, there is only one $B$-stable line in $H^{0}(w,-\lambda-\rho-D)$. Thus, $H^{l(w)}(w,\lambda)$ is a cyclic $B$- module.    
  
So, it is sufficient to prove that $w\cdot \lambda$ is the highest 
weight of $H^{l(w)}(w,\lambda)$. We prove this by induction on $l(w)$.
If $l(w)=0$, there is nothing to prove. 
So, let $w$ be such that $l(w)$ is positive. Let $\phi\in W$ be such that 
$R^{+}(w)\subset R^{+}(\phi)$. Let $\lambda$ be such that $\lambda +\rho$ is 
non singular and $\phi\cdot \lambda$ is dominant. Since $l(w)$ is positive,
there is a simple root $\alpha$ such that $w(\alpha)$ $<$ $0$. Therefore, 
$\phi(\alpha) < 0$. Hence, we have 
$\langle \lambda, \alpha \rangle \leq -2$. Therefore, 
$\langle s_{\alpha}\cdot\lambda, \alpha \rangle \geq 0$.

Now, using the short exact sequence of $B$- modules
$0\longrightarrow K\longrightarrow V_{s_{\alpha}\cdot \lambda, \alpha}
\longrightarrow s_{\alpha}\cdot\lambda\longrightarrow (0)$, we get the  
following exact sequence of $B$- modules:
$H^{l(w)}(w,\lambda)=H^{l(w s_{\alpha})}(w s_{\alpha}, V_{s_{\alpha}\cdot
\lambda, \alpha})= H^{l(w s_{\alpha})}(w, s_{\alpha}\cdot \lambda)
\longrightarrow H^{l(w s_{\alpha})}(w s_{\alpha}, s_{\alpha}\cdot \lambda)
\longrightarrow H^{l(w s_{\alpha})+1}(w s_{\alpha}, K)=(0).$   
Observing that $R^+(w s_{\alpha})\subset R^+(\phi s_{\alpha})$, and 
$l(w s_{\alpha})=l(w)-1$, by induction, it follows that 
the highest weight of the cyclic $B$- module
$H^{l(w s_{\alpha})}(w s_{\alpha}, s_{\alpha}\cdot\lambda)$ is 
$w s_{\alpha}\cdot s_{\alpha}\cdot \lambda=w\cdot \lambda$.

Since the map above $H^{l(w)}(w,\lambda)\longrightarrow 
H^{l(w s_{\alpha})}(w s_{\alpha}, s_{\alpha}\cdot\lambda)$ is surjective, 
the highest weight of $H^{l(w)}(w,\lambda)$ is also $w\cdot \lambda$. 

{\bf Step(2):} We now prove the Proposition by induction on $l(w)$. If $l(w)=0$, 
there is nothing to prove. Assume that the $l(w)$ is positive. So, choose a 
simple root $\gamma$ such that $l(w)=1+l(s_{\gamma}w)$. Then, by the exact 
sequence,$(0)\longrightarrow H^{1}(s_{\gamma}, H^{l(s_{\gamma}w)}(s_{\gamma}w,\lambda))\longrightarrow H^{l(w)}(w,\lambda)\longrightarrow H^{0}(s_{\gamma},
H^{l(w)}(s_{\gamma}w, \lambda))=H^{0}(s_{\gamma},(0))=(0)$, we must have 
$H^{l(s_{\gamma}w)}(s_{\gamma}w,\lambda)\neq (0)$. Therefore, by induction 
$R^{+}(s_{\gamma}w)\subset R^{+}(\phi)$. 
Hence, by the above observation the highest weight of 
the cyclic $B$- module $H^{l(s_{\gamma}w)}(s_{\gamma}w,\lambda)$ is 
$s_{\gamma}w\cdot\lambda$.

Since $H^{1}(s_{\gamma}, H^{l(s_{\gamma}w)}(s_{\gamma}w,\lambda))=H^{l(w)}(w,\lambda)\neq (0)$, the map $H^{l(s_{\gamma}w)}
(s_{\gamma}w,\lambda)\otimes L_{\gamma} \longrightarrow H^{l(w)}(w,\lambda)$ (induced by the evaluation map $H^{0}(s_{\gamma}, 
H^{l(s_{\gamma}w)}(s_{\gamma} w, \lambda)^{*}\otimes L_{-\gamma})\longrightarrow 
H^{l(s_{\gamma}w)}(s_{\gamma} w, \lambda)^{*}\otimes L_{-\gamma}$ and using 
Serre-duality on $\bp^{1}$) is non zero. Therefore,    
if $V$ denote the $B_{\gamma}$- indecomposable component of 
$H^{l(s_{\gamma}w)}(s_{\gamma}w, \lambda)$ containing the 
weight space of weight $s_{\gamma}w\cdot\lambda$, then, 
$H^{1}(s_{\gamma}, V)$ is non zero (since $B$-span of $V$ is 
$H^{l(s_{\gamma}w)}(s_{\gamma}w, \lambda)$). Hence, either 
$\langle s_{\gamma}w\cdot \lambda, \gamma \rangle \leq -1$ or 
$w\cdot\lambda=s_{\gamma}\cdot s_{\gamma}w\cdot\lambda$ is a weight of 
$V$. But, the second possibility is violated from corollary 3.2, since 
$w\not\leq s_{\gamma}w$.Hence, we must have $\langle s_{\gamma}w\cdot\lambda, \gamma\rangle \leq -1$. Now, since $\phi\cdot\lambda=\phi (s_{\gamma}w)^{-1}\cdot s_{\gamma}w\cdot\lambda$ is a dominant weight, and 
$R^{+}(s_{\gamma}w)\subset R^{+}(\phi)$, from the above observation, it 
follows that $\phi(s_{\gamma}w)^{-1}(\gamma)$ is a negative root and hence 
$R^{+}(w)\subset R^{+}(\phi)$. \end{proof}

Let $\phi\in W$ be arbitrary. Let $\lambda$ be a weight such that $\lambda+\rho$ is non singular and $\phi\cdot\lambda$ is dominant. Then, we have 
\bcor The restriction map 
$H^{l(\phi)}(G/B,\lambda)\longrightarrow H^{l(\phi)}(\phi ,\lambda)$
is surjective. \ecor
\begin{proof}  By induction on $l(\phi)$. If $l(\phi)=0$, there is nothing to prove. So, let $l(\phi)$ be a positive integer. Choose a simple root $\alpha$
such that $\phi(\alpha)$ is negative. Hence $\langle \lambda , \alpha
\rangle \leq -2$. Now, using the short exact sequnce of $B$- modules: $(0)\longrightarrow K\longrightarrow V_{s_{\alpha}\cdot\lambda,\alpha}\longrightarrow s_{\alpha}\cdot \lambda \longrightarrow (0)$, we get the following exact sequence of $B$- modules:$H^{l(\phi s_{\alpha})}(\phi s_{\alpha}, V_{s_{\alpha}\cdot\lambda, \alpha})\longrightarrow H^{l(\phi s_{\alpha})}(\phi s_{\alpha}, s_{\alpha}\cdot \lambda) \longrightarrow H^{1+l(\phi s_{\alpha})}(\phi s_{\alpha}, K)=(0)$.
But, we have $H^{l(\phi s_{\alpha})}(\phi s_{\alpha}, V_{s_{\alpha}\cdot 
\lambda, \alpha})=H^{l(\phi s_{\alpha})}(\phi, s_{\alpha}\cdot\lambda)
=H^{l(\phi)}(\phi,\lambda)$. Now, consider the following commutative diagram of $B$- modules: $$\begin{array}{ccc}
H^{l(\phi)}(G/B,\lambda) & \longrightarrow & H^{l(\phi)}(\phi, \lambda) \\
\downarrow              &&                   \downarrow \\
H^{l(\phi s_{\alpha})}(G/B, s_{\alpha}\cdot \lambda) & \longrightarrow 
& H^{l(\phi s_{\alpha})}(\phi s_{\alpha}, s_{\alpha}\cdot\lambda) \\
\end{array}$$ 
The isomorphism of the left vertical map follows from Theorem 8.3.11 in 
\cite{kum2}. The right vertical map is surjective by the above observation.  The second horizontal map is surjective by induction.
Hence, the image of a highest weight vector of weight
$\phi\cdot\lambda$ in $H^{l(\phi)}(G/B, \lambda)$ via the first
horizontal map is non zero.  Now, the surjectivity of the first
horizontal map follows from Proposition 4.2. \end{proof}

Let $\phi\in W$ be arbitrary. Let $M:=max \{\langle \beta,\gamma\rangle :\beta\in \phi(S), \gamma\in S\}$. 
Let $\lambda$ be a weight such that  
$\langle \phi\cdot \lambda, \gamma\rangle > l(\phi)M$ for all $\gamma\in S$. We
first prove  
\bprop Let $w\in W$ be arbitrary. Then, for any integer $i > l(\phi)$, 
we have $H^{i}(w, \lambda)=0$. \eprop 
\begin{proof} Proof is by induction on $l(\phi)$.
If $l(\phi)=0$, then, $\lambda$ is dominant, and so all higher cohomologies
vanish [cf Proposition 8.2.2 in \cite{kum2}] and [cf Lemma 137, Chapter 18 in \cite{Math}]. So, let $\phi\in W$ be such that $l(\phi)$ is positive. 
Since $l(\phi)$ is positive, there is a simple root $\alpha$ such that 
$\phi(\alpha)$ is a negative root. Then, 
$\langle \lambda , \alpha \rangle < 0$. If for such an $\alpha$, $w(\alpha)$ is a negative root, then, $H^{i}(w, \lambda)=H^{i-1}(w, s_{\alpha}\cdot \lambda)=(0)$for $i>l(\phi)$. Otherwise, consider the following exact sequence of $B$ - modules:$(0)\longrightarrow \lambda \longrightarrow V_{s_{\alpha}(\lambda), \alpha}\longrightarrow Q\longrightarrow (0)$.  
This induces the following exact sequence of $B$- modules:
$H^{i-1}(w, Q)\longrightarrow H^{i}(w, \lambda)\longrightarrow H^{i}(ws_{\alpha}, s_{\alpha}(\lambda))$.

Since $l(\phi s_{\alpha})=l(\phi)-1$ and $\langle 
\phi s_{\alpha}\cdot s_{\alpha}(\lambda), \gamma \rangle > l(\phi)(M+\langle
\phi(\alpha), \gamma \rangle)\geq (l(\phi)-1)M$, by induction, we have 
$H^{i}(w, s_{\alpha}(\lambda))=H^{i}(ws_{\alpha}, s_{\alpha}(\lambda))=0$ for 
any $i>l(\phi)-1$.    
Therefore, from the above exact sequence of $B$- modules above, it is 
sufficient to prove that $H^{i-1}(w, Q)=(0)$ for any $i>l(\phi)$. 
To prove this, we consider the following exact sequence of $B$- modules:
$(0)\longrightarrow V_{s_{\alpha}\cdot \lambda, \alpha}\longrightarrow Q
\longrightarrow s_{\alpha}(\lambda)\longrightarrow (0)$. This gives the following exact sequence of $B$- modules:
  
$H^{i-1}(w s_{\alpha}, s_{\alpha}\cdot \lambda)\longrightarrow H^{i-1}(w,Q)
\longrightarrow H^{i-1}(w, s_{\alpha}(\lambda))$.
Here, by induction, $H^{i-1}(w s_{\alpha}, s_{\alpha}\cdot \lambda)=0
=H^{i-1}(w, s_{\alpha}(\lambda))$. Therefore, $H^{i-1}(w, Q)=0$.  
\end{proof}

\section{\bf On the index of the least and the topmost non vanishing cohomology modules}

We first setup some notation to describe the results in this section. We then prove some combinatorial Lemmas. We then prove the following Theorem.
In this section, we prove the following results.
Let $(w, \phi)\in W\times W$ be arbitrary. Let $M=\max\{\langle \beta, \gamma
\rangle ~ : ~ \beta \in \phi(S), ~ \gamma\in S\}$. Let $\tau^{+}(w,\phi)$ (resp.$\tau^{-}(w,\phi)$) be the unique maximal element as in Lemma 5.3(1) (resp. Lemma 5.4(1)). Also, let $D$ denote the boundary divisor of $X(w)$. 
With these notations, we have the following: 

{\bf Theorem}\begin{enumerate}
\item If $\lambda$ is a weight such that $\phi\cdot\lambda$ is dominant,
then, the restriction map $H^{l^{+}(w,\phi)}(w,\lambda)
\longrightarrow H^{l^{+}(w,\phi)}(\tau^{+}(w,\phi),\lambda)$ is surjective 
and in particular $H^{l^{+}(w,\phi)}(w,\lambda)$ is a non zero $B$- module 
with $\tau^{+}(w,\phi)\cdot\lambda$ as a highest weight of this $B$- module. 
\item When $G$ is finite dimensional, and if $\langle \phi\cdot\lambda , \gamma \rangle >  l(w _{0}\phi)M$, for all simple roots $\gamma$, then, we have $H^{i}(w, \lambda)=0$ for $i > l^{+}(w, \phi)$.
\item If $\lambda$ is a weight such that $\langle \phi\cdot\lambda, \gamma \rangle ~ > ~ l(\phi)M$ for all simple roots $\gamma$, then, the cohomologies $H^{i}(w,\lambda)$ vanish for all $i ~ < ~ l(w)-l^{-}(w,\phi)$.
\item If $\lambda$ is a weight such that $\langle \phi(\lambda), \gamma \rangle \geq -1$ for all simple roots, then, $H^{l(w)-l^{-}(w, \phi)}(w,\lambda- D)\neq 0$.\end{enumerate}

\subsection{\bf Relative lengths of $w$ and $\phi$}
We define relative lengths of $w$ and $\phi$.

Notation: For any $w, ~ \phi\in W$, we set a notation 
$W^{+}(w,\phi):=\{\tau\leq w : R^{+}(\tau)\subset R^{+}(\phi)\}$.   
We also set another notation, $W^{-}(w,\phi):=\{\tau\leq w:R^{+}(\tau)\bigcap
R^{+}(\phi) ~ is ~ empty\}$.  
Now, for a given pair $(w,\phi)\in W\times W$, we define two relative lengths 
as follows: 

1. $l^{+}(w,\phi):=\max \{l(\tau):\tau\in W^{+}(w,\phi)\}$.  

2. $l^{-}(w,\phi):=\max \{l(\tau):\tau\in W^{-}(w,\phi)\}$. We have  

\blem Let $\alpha$ be a simple root. Let $\tau,  ~ \phi \in W$ be such that 
both roots $\tau(\alpha)$ and $\phi(\alpha)$ are positive. 
Then, if $R^{+}(\tau)\subset R^{+}(\phi s_{\alpha})$, then,
$R^{+}(\tau)\subset R^{+}(\phi)$. \elem 
\begin{proof} Let $\beta \in R^{+}(\tau)$ be arbitrary. 
Since $\tau(\alpha) >0$, $\beta \neq \alpha$ and so $s_{\alpha}(\beta)>0$.

{\bf Case(1):} If $\tau(s_{\alpha}(\beta))$ is a negative root, then, $\phi(\beta)=\phi s_{\alpha}(s_{\alpha}(\beta))$ is a negative root since $R^{+}(\tau)\subset R^{+}(\phi s_{\alpha})$. Hence, $\beta\in R^{+}(\phi)$. 

{\bf Case(2):} If $\tau(s_{\alpha}(\beta))$ is a positive root, then 
$s_{\alpha}(\beta)= \beta +m\alpha$, with $m$ is a positive integer.
(otherwise, $-m\geq 0$, and so $\tau(s_{\alpha}(\beta))=\tau (\beta) -
m\tau(-\alpha) ~ < ~ 0$, since $\tau(\beta) ~ < ~ 0, ~ \tau(\alpha) ~ > ~ 0$
and $-m\geq 0$.)
Now, $\phi(\beta)=\phi s_{\alpha}(s_{\alpha}(\beta))=\phi s_{\alpha}(\beta)
+m \phi s_{\alpha}(\alpha) ~ < ~ 0$, since $\phi s_{\alpha}(\beta) ~ < ~ 0$,
$\phi(\alpha) ~ > ~ 0$ and $m ~ > ~ 0$. 
Hence, in this case also, $\beta \in R^{+}(\phi)$.
Thus, we have $R^{+}(\tau)\subset R^{+}(\phi)$.\end{proof}

We also have: \blem Let $\alpha$ be a simple root. Let $\tau, ~ \phi\in W$ be such that $\tau(\alpha) > 0$ and $\phi(\alpha)<0$. Then, if $R^{+}(\tau)\bigcap R^{+}(\phi s_{\alpha})$ is empty, then, $R^{+}(\tau)\bigcap R^{+}(\phi)$ is empty. \elem
Proof is similar to that of Lemma(6.1).

The relative length $l^{+}(w, \phi)$ satisfies the following properties:
\blem  
\begin{enumerate}
\item The set $W^{+}(w, \phi)$ has a unique maximal element with respect to the Bruhat order. More precisely, there is a unique $\tau^{+}(w,\phi)\in 
W^{+}(w,\phi)$ such that $l(\tau^{+}(w,\phi))=l^{+}(w,\phi)$ and further, for 
any $\tau\in W^{+}(w,\phi)$,we must have $\tau \leq \tau^{+}(w,\phi)$.
\item For any $\tau \leq w$, and for any $\phi\in W$, we have $l^+(\tau, \phi)\leq l^+(w, \phi)$. 
\item For any $\gamma\in S$ such that $w^{-1}(\gamma) > 0$ and for any $\phi$, $\tau^{+}(s_{\gamma}w,\phi)\in \{\tau^{+}(w,\phi), ~  s_{\gamma}\tau^{+}(w,\phi)\}$ and $l^{+}(w,\phi)\leq l^{+}(s_{\gamma}w, \phi)\leq 1+l^{+}(w, \phi)$.
\item For any simple root $\alpha$, and for any pair $(w, \phi)\in W\times W$ 
such that both the roots $w(\alpha)$ and $\phi(\alpha)$ are positive roots,
the following holds:
\begin{itemize}\item[a.] $\tau^{+}(ws_{\alpha},\phi)=\tau^{+}(w,\phi)$ and  
$l^{+}(ws_{\alpha},\phi) = l^{+}(w,\phi)$.  
\item[b.]$\tau^{+}(ws_{\alpha},\phi s_{\alpha})=\tau^{+}(w,\phi)s_{\alpha}$ and  $l^{+}(w s_{\alpha},\phi s_{\alpha})$ $=$ $1+l^{+}(w,\phi)$. 
\item[c.]$l^{+}(w , \phi s_{\alpha})$ $\leq$ $1+l^{+}(w, \phi)$. 
\end{itemize}
\end{enumerate}\elem
\begin{proof}

We first make the following observation:

{\bf Observation:(1)} If $\alpha$ is a simple root such that $w(\alpha)>0$ and $\phi(\alpha)>0$, then, $l^{+}(w,\phi)\leq l^{+}(ws_{\alpha},\phi s_{\alpha})-1$.
For a proof: If $\tau\in W^{+}(w,\phi)$, then, $\tau(\alpha)>0$. Hence $\tau s_{\alpha}\in W^{+}(ws_{\alpha}, \phi s_{\alpha})$ and $l(\tau)=l(\tau s_{\alpha})-1$. Therefore, $l(\tau)\leq l^{+}(w s_{\alpha},\phi s_{\alpha})-1$. Since $\tau$ was arbitrary in $W^{+}(w,\phi)$, we are done. 

For a proof of $(1)$: We prove this by induction on $l(w)$. If $l(w)=1$, then, $w=s_{\alpha}$ for some simple root $\alpha$. In this case, either $\tau^{+}(w,\phi)=1$ or $\tau^{+}(w,\phi)=s_{\alpha}$ depending on whether $\phi(\alpha) >0$ or $\phi(\alpha) <0$.
So, assume that $l(w)\geq 2$. Choose a simple root $\alpha$ such that 
$w(\alpha)<0$. 
Then, we have two possibilities:

{\bf Case(1):} If $\phi(\alpha) >0$, then, for any $\tau\in W^{+}(w,\phi)$, $\tau(\alpha) >0$ and hence $\tau\leq w s_{\alpha}$, and therefore $\tau\in 
W^{+}(w s_{\alpha}, \phi)$. The other inequality $W^{+}(w s_{\alpha}, \phi)
\subset W^{+}(w, \phi)$ is trivial. Thus, we have $W^{+}(w,\phi)=
W^{+}(w s_{\alpha},\phi)$. Since $l(w s_{\alpha})=l(w)-1$, by induction, 
there is unique $\tau^{+}(w s_{\alpha},\phi)\in W^{+}(w s_{\alpha}, \phi)$
such that $l(\tau^{+}(ws_{\alpha},\phi))=l^{+}(w,\phi)$, for any $\tau\in W^{+}(w s_{\alpha},\phi)$, we must have $\tau\leq \tau^{+}(w s_{\alpha}, \phi)$. Since $W^{+}(w,\phi)= W^{+}(w s_{\alpha}, \phi)$, the assertion $(1)$ is immediate.     

{\bf Case(2):} If $\phi(\alpha) <0$, let $\tau_{0}\in W^{+}(w,\phi)$ be such that $l(\tau_{0})=l^{+}(w,\phi)$. Then, $\tau_{0}(\alpha)<0$, since otherwise, $\tau_{0}s_{\alpha}\in W^{+}(w, \phi)$, by Lemma 5.1, and $l(\tau_{0}s_{\alpha})=l(\tau_{0})+1$, which is a contradiction to the hypothesis that  $l(\tau_{0})=l^{+}(w,\phi)$.

Since $\tau_{0}(\alpha) <0$, $\tau_{0}\in W^{+}(w,\phi)$, 
and both $w(\alpha)$ and $\phi(\alpha)$ are negative roots, it is easy to 
see that $\tau_{0} s_{\alpha}\in W^{+}(w s_{\alpha}, \phi s_{\alpha})$. 
Therefore, using Observation(1), we have $l^{+}(ws_{\alpha},\phi s_{\alpha})\geq l(\tau_{0} s_{\alpha})=l(\tau_{0})-1=l^{+}w,\phi)-1\geq l^{+}(w s_{\alpha}, \phi s_{\alpha})$. Hence, $l^{+}(ws_{\alpha}, \phi s_{\alpha})=l(\tau_{0}s_{\alpha})=l^{+}(w,\phi)-1$.   

On the otherhand, since $l(w s_{\alpha}) = l(w)-1$, by induction,  
there is a unique element $\tau^{+}(w s_{\alpha}, \phi s_{\alpha})\in 
W^{+}(w s_{\alpha}, \phi s_{\alpha})$ such that $l(\tau^{+}(w s_{\alpha},\phi s_{\alpha}))=l^{+}(w s_{\alpha}, \phi s_{\alpha})$, and further it is the unique maximal element with respect to the Bruhat order. Hence, by the uniqueness of the element in $W^{+}(w s_{\alpha}, \phi s_{\alpha})$ having length equal to $l^{+}( w s_{\alpha}, \phi s_{\alpha})$ and from the above observation that $\tau_{0}s_{\alpha}\in W^{+}(ws_{\alpha},\phi s_{\alpha})$, we must have $\tau_{0} s_{\alpha} = \tau^{+}(w s_{\alpha}, \phi s_{\alpha})$. 

We now claim that $\tau_{0}$ is the unique maximal element in $W^{+}(w,\phi)$
with respect to the Bruhat order. Now, for any $\tau\in W^{+}(w,\phi)$, 
we have two possibilities: 

{\bf Subcase(1):} If $\tau(\alpha)$ is a positive root, then, by Lemma 5.1,
$\tau\in W^{+}(w s_{\alpha}, \phi s_{\alpha})$. Now, since 
$\tau_{0} s_{\alpha}$ is the unique maximal element of 
$W^{+}(w s_{\alpha}, \phi s_{\alpha})$ with respect to the Bruhat order, we must have $\tau\leq \tau_{0} s_{\alpha}\leq \tau_{0}$. Thus, we are done in this case.

{\bf Subcase(2):} If $\tau(\alpha)$ is a negative root, then, it is easy to see 
that $\tau s_{\alpha}\in W^{+}(w s_{\alpha}, \phi s_{\alpha})$, and hence
$\tau s_{\alpha}\leq \tau_{0} s_{\alpha}$. Now, since 
$\tau_{0} s_{\alpha}\leq \tau_{0}$, we must have $\tau \leq \tau_{0}$.
Thus, we are done.  

Proof of (2): This follows from the definition of $l^{+}(w,\phi)$. 

Proof of (3): Let 
$\tau \leq s_{\gamma}w$ and $R^{+}(\tau)\subset R^{+}(\phi)$. 
Now, if $\tau^{-1}(\gamma) < 0$, then, $s_{\gamma}\tau\leq w$, and it 
satisfies $R^{+}(s_{\gamma}\tau)\subset R^{+}(\phi)$. Otherwise, 
$\tau\leq w$. The proof now follows from $(2)$. 

Proof of $4(a)$: From the proof of $(1)$, it is easy to see that $W^{+}(w,\phi)=W^{+}(w s_{\alpha},\phi)$ if $\phi(\alpha)> 0$. The assertion $4(a)$ is immediate from this observation. 

Proof of $4(b)$: From the proof of $(1)$, it is easy to see that $1+l^{+}(w,\phi)=l^{+}(w s_{\alpha},\phi s_{\alpha})$. Also, it is easy to see that $\tau^{+}(w,\phi)(\alpha)>0$, and $\tau^{+}(w,\phi)s_{\alpha})\in W^{+}(w,\phi)$. 
Hence, we have $\tau^{+}(w,\phi)s_{\alpha}=\tau^{+}(ws_{\alpha},\phi s_{\alpha})$.

Proof of $4(c)$: This follows from $4(b)$ and $(2)$. \end{proof}

We also have the following Lemma:
\blem 
\begin{enumerate}
\item The set $W^{-}(w,\phi)$ has a unique maximal element with respect to the Bruhat order. More precisely, there is a unique $\tau^{-}(w,\phi)\in W^{-}(w,\phi)$ such that $l(\tau^{-}(w,\phi))=l^{-}(w,\phi)$ and further, for any $\tau\in 
W^{-}(w,\phi)$, we must have $\tau\leq \tau^{-}(w,\phi)$.
\item For any $\tau\leq w$, and for any $\phi\in W$, we have 
$l^{-}(\tau, \phi)\leq l^{-}(w,\phi)$.
\item If $\gamma\in S$ is such that $w^{-1}(\gamma)>0$, then, 
for any $\phi\in W$, we have $\tau^{-}(s_{\gamma}w,\phi)\in \{\tau^{-}(w,\phi), s_{\gamma}\tau^{-}(w,\phi)\}$ and $l^{-}(w, \phi)\leq l^{-}(s_{\gamma}w, \phi)\leq 1+l^{-}(w,\phi)$. 
\item For any simple root $\alpha$, and for any pair $(w,\phi)\in W\times W$
such that $w(\alpha) > 0$ and $\phi(\alpha) < 0$, then, the following holds:
\begin{itemize} 
\item[a.] $\tau^{-}(w,\phi)=\tau^{-}(ws_{\alpha},\phi)$ and $l^{-}(w s_{\alpha}, \phi)=l^{-}(w,\phi)$.
\item[b.] $\tau^{-}(ws_{\alpha},\phi s_{\alpha})=\tau^{-}(w,\phi)s_{\alpha}$ and $l^{-}(w s_{\alpha}, \phi s_{\alpha})=1+l^{-}(w,\phi)$.
\item[c.] $l^{-}(w ,\phi s_{\alpha})\leq 1+l^{-}(w,\phi)$.
\end{itemize}
\end{enumerate}\elem

Proof of the Lemma is similar to that of Lemma 5.3.

Let $(w,\phi)\in W\times W$ be arbitrary. We then have:

\blem $\tau^{+}(w,\phi)\phi^{-1}\leq \tau\phi^{-1}$ for any $\tau\leq w$. \elem
\begin{proof} {\bf Step(1):} We first show that $\tau^{+}(w,\phi)\phi^{-1}\leq w\phi^{-1}$. 

{Proof of Step(1):} Proof is by induction on $l(w)$.

If $l(w)=0$, the assertion is trivial.
So, let $l(w)\geq 1$. Then, there exists a simple root $\alpha$ such that $w(\alpha) < 0$. Now, we have two possibilities:

{\bf Case(1):} $\phi(\alpha)<0$. Then, by Lemma(5.3(4(b)), we have $\tau^{+}(w,\phi)\phi^{-1}=\tau^{+}(ws_{\alpha}, \phi s_{\alpha})(\phi s_{\alpha})^{-1}$. Since $l(ws_{\alpha})=l(w)-1$, by induction, $\tau^{+}(ws_{\alpha}, \phi s_{\alpha})(\phi s_{\alpha})^{-1}\leq ws_{\alpha}(\phi s_{\alpha})^{-1}=w\phi^{-1}$. Thus, from above observation , we have $\tau^{+}(w,\phi)\phi^{-1}\leq w\phi^{-1}$.

{\bf Case(2):} $\phi(\alpha)>0$. Then, by Lemma(6.3(4(a)), we have $\tau^{+}(w,\phi)\phi^{-1}=\tau^{+}(ws_{\alpha},\phi)\phi^{-1}$. Since $l(ws_{\alpha})=l(w)-1$, by induction, $\tau^{+}(w,\phi)\phi^{-1}=\tau^{+}(w s_{\alpha},\phi)\phi^{-1}\leq ws_{\alpha}\phi^{-1}$. But, on the other hand, since $\phi(\alpha)$ is a real positive root with $\phi(\alpha)\in R^{+}(w\phi^{-1})$, we must have $ws_{\alpha}\phi^{-1}=w\phi^{-1}s_{\phi(\alpha)} < w\phi^{-1}$. Hence, we have $\tau^{+}(w,\phi)\phi^{-1}\leq ws_{\alpha}\phi^{-1} < w\phi^{-1}$. This completes the proof of Step(1).

{\bf Step(2):} We show that for any $w\in W$, $\tau\in W^{+}(w,\phi)$, $\tau \phi^{-1}\geq \tau^{+}(w,\phi)\phi^{-1}$.

{Proof of Step(2):} Proof is by induction on $l(w)$. 

If $l(w)=0$, there is nothing to prove. So, let $l(w)\geq 1$. Then, there is a simple root such that $w(\alpha) < 0$. Then, we have two possibilities:

{\bf Case(1):} $\phi(\alpha) < 0$. Let $\tau\in W^{+}(w,\phi)$. Then, we have two possibilities.

{\bf Subcase(i):} $\tau(\alpha) > 0$. Then, we must have $\tau\in W^{+}(w s_{\alpha}, \phi s_{\alpha})$. Since $l(w s_{\alpha})=l(w)-1$, by induction and using the Lemma(5.3(4(b))), we have $\tau (\phi s_{\alpha})^{-1}\geq \tau^{+}(ws_{\alpha}, \phi s_{\alpha})(\phi s_{\alpha})^{-1}=\tau^{+}(w,\phi)\phi^{-1}$. But, on the other hand, since $-\phi(\alpha)\in R^{+}(\tau \phi^{-1})$, we have $\tau (\phi s_{\alpha})^{-1}=\tau \phi^{-1}s_{-\phi(\alpha)} < \tau \phi^{-1}$. Thus, we have $\tau \phi^{-1} >  \tau (\phi s_{\alpha})^{-1} \geq \tau^{+}(w,\phi) \phi^{-1}$.
  
{\bf Subcase(ii):} $\tau(\alpha) < 0$. Then, we must have $\tau s_{\alpha}\in W^{+}(w s_{\alpha}, \phi s_{\alpha})$. Since $l(w s_{\alpha})=l(w)-1$, by induction and using Lemma(5.3(4(b))), we have $\tau \phi^{-1}=\tau s_{\alpha}(\phi s_{\alpha})^{-1}\geq \tau^{+}(w s_{\alpha}, \phi s_{\alpha})(\phi s_{\alpha})^{-1}=\tau^{+}(w,\phi)\phi^{-1}$. 

{\bf Case(2):} $\phi(\alpha) > 0$. By Lemma(5.3(4(a))), we have $\tau^{+}(w,\phi)=\tau^{+}(w s_{\alpha}, \phi)$ and $W^{+}(w,\phi)=W^{+}(w s_{\alpha},\phi)$. Now, let $\tau\in W^{+}(w,\phi)=W^{+}(w s_{\alpha},\phi)$. Since $l(w s_{\alpha})=l(w)-1$, and since $\tau\in W^{+}(w s_{\alpha},\phi)$, by induction, we must have $\tau\phi^{-1}\geq\tau^{+}(w s_{\alpha},\phi)\phi^{-1}=\tau^{+}(w,\phi)\phi^{-1}$. This completes the proof of Step(2).

{\bf Step(3):} We now prove the Lemma. Let $\tau\leq w$. Then, $\tau^{+}(\tau,\phi)\in W^{+}(w,\phi)$. Hence, by Step(2), we must have $\tau^{+}(w,\phi)\phi^{-1}\leq \tau^{+}(\tau,\phi)\phi^{-1}$.
Also,  by Step(1), we have $\tau^{+}(\tau,\phi)\phi^{-1}\leq\tau\phi^{-1}$.   
Thus, we have $\tau^{+}(w,\phi)\phi^{-1}\leq \tau\phi^{-1}$. 
\end{proof}

We next state with out proof a similar Lemma for $\tau^{-}(w,\phi)\phi^{-1}$.

\blem For any $\tau\leq w$, $\tau^{-}(w,\phi)\phi^{-1}\geq \tau\phi^{-1}$ with respect to the Bruhat order.\elem

{\bf Proof} Similar to  that of Lemma 5.5.

\bcor
Let $(w,\phi)\in W\times W$. Let $\lambda$ be a weight such that $\phi\cdot\lambda$ is dominant. Let $i$ be a non-negative integer. Then, any weight $\mu$ of $H^{i}(w,\lambda)$ satisfies $\tau^{-}(w,\phi)\cdot\lambda\leq \mu\leq \tau^{+}(w,\phi)\cdot\lambda$ \ecor
\begin{proof} For any $\tau\leq w$ from Lemma 5.5, we have $\tau^{+}(w,\phi)\phi^{-1}\leq \tau\phi^{-1}$. Since $\phi\cdot\lambda$ is dominant, for any $\tau\leq w$, $\tau\cdot\lambda=\tau\phi^{-1}\cdot\phi\cdot\lambda\leq \tau^{+}(w,\phi)\phi^{-1}\cdot\phi\cdot\lambda=\tau^{+}(w,\phi)\cdot\lambda$. Similarly, for any $\tau\leq w$, using Lemma 5.6, we have $\tau\cdot\lambda\geq\tau^{-}(w,\phi)\cdot\lambda$. Thus, we have $\tau^{-}(w,\phi)\cdot\lambda\leq\tau\cdot\lambda\leq\tau^{+}(w,\phi)\cdot\lambda$ for any $\tau\leq w$.      

Since any  weight $\mu$ of $H^{i}(w,\lambda)$ is in the convex hull of the set $\{\tau\cdot\lambda: \tau\leq w\}$, the assertion follows from the above observation. \end{proof}

Let $(w, \phi)\in W\times W$ be arbitrary. Let $M=\max\{\langle \beta, \gamma
\rangle ~ : ~ \beta \in \phi(S), ~ \gamma\in S\}$. Let $\tau^{+}(w,\phi)$ (resp.$\tau^{-}(w,\phi)$) be the unique maximal element as in Lemma 5.3(1) (resp. Lemma 5.4(1)). Also, let $D$ denote the boundary divisor of $X(w)$. 
With these notations, we have the following: 
\bth \begin{enumerate}
\item If $\lambda$ is a weight such that $\phi\cdot\lambda$ is dominant,
then, the restriction map $H^{l^{+}(w,\phi)}(w,\lambda)
\longrightarrow H^{l^{+}(w,\phi)}(\tau^{+}(w,\phi),\lambda)$ is surjective 
and in particular $H^{l^{+}(w,\phi)}(w,\lambda)$ is a non zero $B$- module 
with $\tau^{+}(w,\phi)\cdot\lambda$ as a highest weight of this $B$- module. 
\item When $G$ is finite dimensional, and if $\langle \phi\cdot\lambda , \gamma \rangle >  l(w _{0}\phi)M$, for all simple roots $\gamma$, then, we have $H^{i}(w, \lambda)=0$ for $i > l^{+}(w, \phi)$.
\item If $\lambda$ is a weight such that $\langle \phi\cdot\lambda, \gamma \rangle ~ > ~ l(\phi)M$ for all simple roots $\gamma$, then, the cohomologies $H^{i}(w,\lambda)$ vanish for all $i ~ < ~ l(w)-l^{-}(w,\phi)$.
\item If $\lambda$ is a weight such that $\langle \phi(\lambda), \gamma \rangle \geq -1$ for all simple roots, then, $H^{l(w)-l^{-}(w, \phi)}(w,\lambda- D)\neq 0$.
\end{enumerate}
\eeth
\begin{proof}{\bf Proof of (1):}
We now fix $\phi$ but arbitrary. Let $\lambda$ be a weight as in the 
hypothesis of the Theorem. By induction on $l(\phi)-l^{+}(w,\phi)$, we prove 
that the restriction map $H^{l^{+}(w,\phi)}(w,\lambda)\longrightarrow 
H^{l^{+}(w,\phi)}(\tau^{+}(w,\phi) ,\lambda)$ is onto. Then, from Proposition 4.2, and from the fact that $l^{+}(w,\phi)=l(\tau^{+}(w,\phi))$, $H^{l^{+}(w,\phi)}(w,\lambda)$ is non zero and $\tau^{+}(w,\phi)\cdot\lambda$ is a weight of it.
From Corollary 5.7, $\tau^{+}(w,\phi)\cdot\lambda$ is a highest weight of 
$H^{l^{+}(w,\phi)}(w,\lambda)$.

If $l(\phi)-l^{+}(w,\phi)=0$, then $\phi \leq w$. In this case, the assertion 
follows from Corollary 4.3. 
So, let $w$ be such that $l(\phi)-l^{+}(w, \phi)$ is a positive integer.
Consider $\tau^{+}(w,\phi)$. Since $l(\phi)-l^{+}(w,\phi)$ is a positive 
integer, we must have $l(\tau^{+}(w,\phi))=l^{+}(w,\phi)\leq l(\phi)-1$.
Therefore, there is a simple root $\gamma$ such that $s_{\gamma}
\tau^{+}(w,\phi) > \tau^{+}(w,\phi)$
and $R^{+}(s_{\gamma}\tau^{+}(w,\phi))\subset R^{+}(\phi)$. 
By maximality of $\tau^{+}(w,\phi)$, $s_{\gamma}w > w$.(otherwise, 
$s_{\gamma}\tau^{+}(w,\phi)\leq w$, contradicting the
maximality of $\tau^{+}(w,\phi)$ with this property.)
Hence, from the Lemma 5.3(3) and from the fact that $l(s_{\gamma}\tau^{+}(w,\phi))=1+l(\tau^{+}(w,\phi))$, we must have $l^{+}(s_{\gamma}w, \phi)=l(s_{\gamma}
\tau^{+}(w,\phi)) =1+l(\tau^{+}(w,\phi))
=1+l^{+}(w, \phi)$, and therefore $s_{\gamma}\tau^{+}(w,\phi)=
\tau^{+}(s_{\gamma}w,\phi)$. 

Let $V_{1}:= H^{1}(s_{\gamma}, H^{l^{+}(w, \phi)}(w, \lambda))$, 
$V:= H^{1+l^{+}(w, \phi)}(s_{\gamma}w, \lambda)$, 
$V_{2}:= H^{0}(s_{\gamma}, H^{1+l^{+}(w, \phi)}(w, \lambda)$, 
$V_{1}^{\prime}:=H^{1}(s_{\gamma}, H^{l^{+}(w,\phi)}(\tau^{+}(w,\phi),
\lambda)$, $V^{\prime}:=H^{l^{+}(s_{\gamma}w,\phi)}(\tau^{+}(s_{\gamma}w,\phi), \lambda)$, and $V_{2}^{\prime}:=H^{0}(s_{\gamma}, H^{l^{+}(s_{\gamma}w, \phi)}(\tau^{+}(w,\phi),\lambda))$, we have 
the following commutative diagram of $B$- modules:
$$\begin{array}{ccccccccc}
(0) & \rightarrow & V_{1} & \rightarrow & V & \rightarrow 
& V_{2} & \rightarrow & (0)\\  
\downarrow & & \downarrow & & \downarrow & & \downarrow & & \downarrow  \\
(0) & \rightarrow & V_{1}^{\prime}& \rightarrow & V^{\prime}& \rightarrow & 
V_{2}^{\prime}& \rightarrow & (0) \\  
\end{array}.$$
From the observation made above , we have $s_{\gamma} \tau^{+}(w,\phi)
=\tau^{+}(s_{\gamma}w,\phi)$ and  $l(\phi)-l^{+}(s_{\gamma}w,\phi) < 
l(\phi)-l^{+}(w,\phi)$, by induction, the middle vertical arrow 
$V\longrightarrow V^{\prime}$ is surjective. Since $l^{+}(s_{\gamma}w,\phi)=l(\tau^{+}(s_{\gamma}w,\phi))$, from Proposition 4.2, $V^{\prime}=H^{l^{+}(s_{\gamma}w,\phi)}(\tau^{+}(s_{\gamma}w,\phi),\lambda)$ is a cyclic $B$- module with highest weight $\tau^{+}(s_{\gamma}w,\phi)\cdot\lambda$.

We now claim that a highest weight vector $v^{+}$ of weight 
$\tau^{+}(s_{\gamma}w,\phi)\cdot\lambda$ lies in the image of the first 
vertical map $V_{1}\longrightarrow V_{1}^{\prime}\hookrightarrow V^{\prime}$.
Since the restriction map $V\longrightarrow V^{\prime}$ is surjective,
there is a weight vector $v$ of weight $\tau^{+}(s_{\gamma}w,\phi)\cdot\lambda$
in $V$. Now, such a weight vector $v$ must lie in the kernal of the 
horizontal map $V\longrightarrow V_{2}$ in the first row. 
(For a proof of this fact: If the image of $v$ is non zero, then, by Verma Lemma[cf \cite{bwbk}], there is an indecomposbale $B_{\gamma}$- component $V_{v}$ in $H^{1+l^{+}(w,\phi)}(w,\lambda)$ such that $v$ is in 
$H^{0}(s_{\gamma}, V_{v})$. Since the weight of $v$ is 
$\tau^{+}(s_{\gamma}w,\phi)\cdot\lambda$, the highest weight of the cyclic 
$B_{\gamma}$- module $V_{v}$ is 
of the form $\tau^{+}(s_{\gamma}w, \phi)\cdot\lambda+m\gamma$ for some non 
negative integer $m$. Therefore, $\tau^{+}(s_{\gamma}w,\phi)\cdot\lambda+m\gamma$ is a weight of $H^{1+l^{+}(w,\phi)}(w,\lambda)$. Since 
$\phi(\tau^{+}(s_{\gamma}w,\phi))^{-1}(\gamma)$ is a positive root, by the 
Corollary 3.3, $\tau^{+}(s_{\gamma}w,\phi)\leq w$, a contradiction to the observation above that $\tau^{+}(s_{\gamma}w,\phi)=s_{\gamma}\tau^{+}(w,\phi)\not\leq w$.) Since the image of a highest weight vector $v$ in $V_{2}$ is zero, $v$ must lie in the image of 
$V_{1}\longrightarrow V_{1}^{\prime}\hookrightarrow V^{\prime}$. 

We now claim that a highest weight vector $v_{1}$ of weight 
$\tau^{+}(w,\phi)\cdot\lambda$ in $H^{l^{+}(w,\phi)}(\tau^{+}(w,\phi), \lambda)$ lies in the image of the restriction map $H^{l^{+}(w,\phi)}(w,\lambda)\longrightarrow  H^{l^{+}(w,\phi)}(\tau^{+}(w,\phi),\lambda)$, which inparticular 
implies that the restriction is surjective.

For a proof of this claim: Consider the $B$- module $H^{l^{+}(w,\phi)}(\tau^{+}(w,\phi),\lambda)$.Inparticular, this is a $B_{\gamma}$- module.  Decompose this into indecomposable $B_{\gamma}$- submodules. Since 
$V_{1}^{\prime}=H^{1}(s_{\gamma}, H^{l^{+}(w,\phi)}(\tau^{+}(w,\phi),\lambda))$ is non zero, the natural map $H^{l^{+}(w,\phi)}(\tau^{+}(w,\phi),\lambda)\otimes L_{\gamma}\longrightarrow H^{1}(s_{\gamma}, 
H^{l^{+}(w,\phi)}(\tau^{+}(w,\phi), \lambda))=V_{1}^{\prime}$ is non zero.
Since $H^{l^{+}(w,\phi)}(\tau^{+}(w,\phi), \lambda)$ is a cyclic $B$- module with highest weight $\tau^{+}(w,\phi)\cdot\lambda$, there is an $B_{\gamma}$-
indecomposable component $U_{1}$ of 
$H^{l^{+}(w,\phi)}(\tau^{+}(w,\phi),\lambda)$ containing a weight vector of weight $\tau^{+}(w,\phi)\cdot\lambda$ such that 
$H^{1}(s_{\gamma}, U_{1})$ is non zero. Hence, any weight in 
$H^{1}(s_{\gamma}, U_{1})$ is in the convex hull of 
$\tau^{+}(w,\phi)\cdot\lambda-r\gamma$ and $s_{\gamma}\cdot(\tau^{+}(w,\phi)\cdot\lambda-r\gamma)$, where $\tau^{+}(w,\phi)\cdot\lambda-r\gamma$ is the lowest 
weight of the cyclic $B_{\gamma}$- module $U_{1}$, with $r$ a non negative 
integer. But, since $V_{1}^{\prime}$ being a highest weight $B$- module with 
highest weight $s_{\gamma}\cdot\tau^{+}(w,\phi)\cdot\lambda$, $r$ must be zero,
and so $U_{1}$ is just the one dimensional vector space $\bc\cdot v_{1}$ 
spanned by $v_{1}$. $H^{1}(s_{\gamma}, \bc\cdot v_{1})$ is a direct summand of 
$V_{1}^{\prime}$ as a $B_{\gamma}$- module and it contains $v^{+}$. 

Now, cosider the commutative diagram of $B_{\gamma}$- modules:      
 $$\begin{array}{ccccc}
H^{l^{+}(w,\phi)}(w,\lambda)\otimes L_{\gamma} & \longrightarrow & 
H^{l^{+}(w,\phi)}(\tau^{+}(w,\phi),\lambda)\otimes L_{\gamma}\\  
\downarrow & & \downarrow \\
V_{1} & \longrightarrow & V_{1}^{\prime}\\  
\end{array}.$$
Since a highest weight vector $v^{+}$ of weight 
$s_{\gamma}\cdot\tau^{+}(w,\phi)\cdot\lambda$ lies in the image of the second
horizontal map followed by the left vertical map, and it also lies in  
the $B_{\gamma}$- direct summand $H^{1}(s_{\gamma}, \bc\cdot v_{1})$, 
$v_{1}\otimes 1$ must lie in the image of the first horizontal map.    
Thus, the restriction map $H^{l^{+}(w,\phi)}(w,\lambda)\longrightarrow 
H^{l^{+}(w,\phi)}(\tau^{+}(w,\phi),\lambda)$ is surjective. 

{\bf Proof of (2):} Here $mathcal{G}$ is finite dimensional and $w_{0}$ denote the 
longest element of the Weyl group.
Proof of $(2)$ is by induction on the length of $w_{0}\phi$. 
If $l(w_{0}\phi)=0$, then, $\tau^{+}(w,\phi)=w$ and so $l^{+}(w,\phi)=l(w)$,
and in this case the assertion follows from the fact that any cohomolgy 
vanishes beyond dimension. 
So, let us assume that $l(w_{0}\phi)$ is positive.
Since $l(w_{0}\phi)$ is positive,  there is a simple root $\alpha$ such 
that $\phi(\alpha)$ is a positive root.
We have two possibilities:

{\bf Case(1):} Now, for such an $\alpha$, if $w(\alpha)$ is a negative root, then, we have $H^{i}(w, \lambda)=H^{i+1}(w, s_{\alpha}\cdot\lambda)$, since $\phi(\alpha)$ is positive. Also, for any $\gamma\in S$, the weight $s_{\alpha}\cdot\lambda$ satisfies $$\langle \phi s_{\alpha}\cdot s_{\alpha}\cdot \lambda, 
\gamma \rangle = \langle \phi\cdot \lambda, \gamma \rangle > l(w_{0}\phi)M
=(1+l(w_{0}\phi s_{\alpha}))M ~ > ~ l(w_{0}\phi s_{\alpha})M.$$ Hence, by 
induction on $\phi$, $H^{i+1}(w, s_{\alpha}\cdot\lambda)=0$ for 
$i+1 > l^{+}(w, \phi s_{\alpha})$. But, from Lemma 5.3(4(a) and 5.3(4(b)), we have $1+l^{+}(w, \phi)= l^{+}(w, \phi s_{\alpha})$. Now, since $i ~ > ~ l^{+}(w,\phi)$, $i+1 ~ > ~ l^{+}(w, \phi s_{\alpha})$ and so by induction on 
$l(w_{0}\phi)$, we must have $H^{i}(w, \lambda)=0$ for all $i > l^{+}(w,\phi)$.

{\bf Case(2):} If for such an $\alpha$, $w(\alpha)$ is a positive root,
then, using the Demazure exact sequence, we get a long exact sequence
$$H^{i}(ws_{\alpha}, \lambda)\longrightarrow H^{i}(w, \lambda)\longrightarrow 
H^{i+1}(w, K)\longrightarrow H^{i+1}(ws_{\alpha}, \lambda).$$ 
Now, by Lemma 5.3(4(a)), we have $l^{+}(w,\phi)=l^{+}(ws_{\alpha},\phi)$ since 
$\phi(\alpha) > 0$. Therefore, by induction on $l(w_{0}w)$, and from the fact that $l(w_{0}ws_{\alpha})=l(w_{0}w)-1=r$, we must have 
$H^{i}(ws_{\alpha}, \lambda)=H^{i+1}(w s_{\alpha}, s_{\alpha}\cdot\lambda)=0$ 
for $i > l^{+}(w, \phi)=l^{+}(ws_{\alpha}, \phi)$, and 
$H^{i+1}(ws_{\alpha}, \lambda)=H^{i+2}(w s_{\alpha}, s_{\alpha}\cdot\lambda)
=0$ for $i > l^{+}(w, \phi)=l^{+}(ws_{\alpha}, \phi)$. 
Hence, from the above exact sequence, we have an isomorphism:
$$H^{i}(w,\lambda)\simeq H^{i+1}(w, K) ~ for ~ i > l^{+}(w, \phi).$$ 
Therefore, to prove the assertion, it is sufficient to prove that 
$H^{i+1}(w, K)=0$ for $i > l^{+}(w, \phi)$. To prove this statement,
we use the other Demazure exact sequence and get the exact sequence :
$$H^{i+1}(w, s_{\alpha}\lambda)\longrightarrow H^{i+1}(w, K)\longrightarrow 
H^{i+1}(ws_{\alpha},\lambda-\alpha)=H^{i+2}(ws_{\alpha}, s_{\alpha}(\lambda)).
$$
Now, by the hypothesis of $\lambda$, for any $\gamma\in S$, 
$s_{\alpha}(\lambda)$ must satisfy 
$$\langle \phi s_{\alpha}\cdot s_{\alpha}(\lambda) , \gamma \rangle 
= \langle \phi\cdot \lambda - \phi(\alpha), \gamma \rangle  > l(w_{0}\phi)M -
\langle \phi(\alpha), \gamma \rangle \geq (l(w_{0}\phi)-1)M=l(w_{0}\phi s_{\alpha})M.$$
Also, by Lemma 5.3(4(a)), we have $l^{+}(w, \phi s_{\alpha})\leq l^{+}(ws_{\alpha}, \phi s_{\alpha})=1+l^{+}(w, \phi)$. since $l(w_{0}\phi s_{\alpha})=l(w_{0}\phi)-1$, by induction , we must have 
$H^{i+1}(w, s_{\alpha}(\lambda))=(0)$ and $H^{i+2}(ws_{\alpha}, s_{\alpha}
(\lambda))=0$ for all $i > l^{+}(w, \phi)$.
Thus, $H^{i+1}(w, K)$ is zero.

{\bf Proof of (3):} The genericity assumption on $\lambda$ in the proof of $(3)$ is $\langle \phi\cdot\lambda, \gamma \rangle ~ > ~ l(\phi)M$.
We now prove $(3)$ by induction on the length of $\phi$.
If $l(\phi)=0$, then $\phi=1$, and $H^{i}(w,\lambda)$ is 
zero for $i < l(w)-l^{-}(w,1)=0$. So, let $\phi$ be such that 
$l(\phi)>0$ and let $\lambda$ be such that $\langle \phi\cdot \lambda, \gamma 
\rangle ~ > ~ l(\phi)M$. 
Since $l(\phi)~ > ~ 0$, there is a simple root $\alpha$ such that 
$\phi(\alpha)$ is a negative root. 

We now prove the assertion by considering two cases.

{\bf Case(1)}: If $w(\alpha)$ is a negative root, then, $H^{i}(w, \lambda)=H^{i-1}(w, s_{\alpha}\cdot \lambda)$ for all $i$. 
We now observe that $H^{i-1}(w,s_{\alpha}\cdot\lambda)=0$ for $i < l(w)-
l^{-}(w, \phi)$, which therefore proves the assertion.
To do this, if $i ~ < ~ l(w)-l^{-}(w, \phi)$, then $i-1 ~ < ~ l(w)- 
l^{-}(w,\phi s_{\alpha})= l(w)-(1+l^{-}(w,\phi))$, by Lemma 5.4(4). Since 
$s_{\alpha}\cdot\lambda$ satisfies $\langle \phi s_{\alpha}\cdot s_{\alpha}\cdot \lambda, \gamma \rangle =\langle\phi\cdot \lambda , \gamma \rangle ~ > ~ l(\phi)M=(l(\phi s_{\alpha})+1)M ~ > ~ l(\phi s_{\alpha})M$, by induction on the 
$l(\phi)$, we have $H^{j}(w,s_{\alpha}\cdot\lambda)=0$ for 
$j < l(w)-l^{-}(w,\phi)$.

{\bf Case(2)}: If $w(\alpha)$ is a positive root, then, consider the following 
exact sequence of $B$ - modules:
$$0\longrightarrow \lambda\longrightarrow V_{s_{\alpha}(\lambda),\alpha}\longrightarrow Q\longrightarrow 0. $$ 
Using this short exact sequence, we get a long exact sequence 
$H^{i-1}(w, Q)\longrightarrow H^{i}(w, \lambda)\longrightarrow 
H^{i}(ws_{\alpha}, s_{\alpha}(\lambda))$. Here, by Lemma 5.4(4), we have 
$l(w)-l^{-}(w,\phi)=l(w s_{\alpha})-l^{-}(ws_{\alpha},\phi s_{\alpha})$, and 
hence, for any $i < l(w)-l^{-}(w, \phi)$, by induction, we must have  
$H^{i}(ws_{\alpha}, s_{\alpha}(\lambda))=0$ since $s_{\alpha}(\lambda)$ 
satisfies $\langle \phi s_{\alpha}\cdot s_{\alpha}(\lambda), \gamma \rangle > l(\phi s_{\alpha}) M$.
So, to prove the assertion, from the above long exact sequence,
it is sufficient to prove that $H^{i-1}(w, Q)$ is zero if  
$i < l(w)-l^{-}(w, \phi)$. To do this, 
we consider the short exact sequence:
$0\longrightarrow V_{s_{\alpha}(\lambda)-\alpha, \alpha}\longrightarrow 
Q\longrightarrow s_{\alpha}(\lambda)\longrightarrow 0$.   

From this short exact sequence, we get a long exact sequence of $B$- modules 
$$H^{i-1}(ws_{\alpha},s_{\alpha}\cdot \lambda)\longrightarrow H^{i-1}(w, Q)
\longrightarrow H^{i-1}(w, s_{\alpha}(\lambda)).$$

Since both the weights $s_{\alpha}\cdot \lambda$ and $s_{\alpha}(\lambda)$ 
satisfy $\langle \phi s_{\alpha}\cdot s_{\alpha}\cdot \lambda, \gamma \rangle 
~ > ~ l(\phi s_{\alpha})M$ and $\langle \phi s_{\alpha}\cdot 
s_{\alpha}(\lambda), \gamma \rangle ~ > ~ l(\phi s_{\alpha})M$ for all 
simple roots $\gamma$, and $l^{-}(w, \phi s_{\alpha})\leq
l^{-}(w s_{\alpha},\phi s_{\alpha})=1+l^{-}(w,\phi)$ (from Lemma 5.4(4)), for  
any $i ~ < ~ l(w)-l^{-}(w, \phi)$, we must have $i-1 ~ < ~ l(w)-
l^{-}(ws_{\alpha}, \phi s_{\alpha})\leq l(w)-l^{-}(w, \phi s_{\alpha})$.

Hence, by induction, for any $i ~ < ~ l(w) - l^{-}(w, \phi)$, we have 
$$H^{i-1}(w s_{\alpha}, s_{\alpha}\cdot \lambda)
=H^{i-1}(w, s_{\alpha}(\lambda))=0.$$ Therefore, the vanishing of 
$H^{i-1}(w, Q)$ follows from the above exact sequence. 

{\bf Proof of $(4)$:}  

{\bf Step(1):} We first show that if $\mu$ is a weight such that $\phi\cdot\mu$ is  negative dominant , then, the restriction map    
$H^{l^{-}(w,\phi)}(w,\lambda)\longrightarrow 
H^{l^{-}(w,\phi)}(\tau^{-}(w,\phi),\mu)$ is surjective by induction on 
$l(\phi)$. 
If $l(\phi)=0$, then $\tau^{-}(w,\phi)=w$ itself and so there is nothing 
to prove. So, let $l(\phi)$ be positive. Choose a simple root $\alpha$ such 
that $\phi(\alpha)$ is a negative root. For such a simple root $\alpha$, we 
have two possibilities:

{\bf Case(1):} $w(\alpha)$ is a negtive root.  
Observation(1):
Since $l(\phi s_{\alpha}) = l(\phi)-1$, by induction, the restriction  
$H^{l^{-}(w,\phi s_{\alpha})}(w,s_{\alpha}\cdot\mu)\longrightarrow 
H^{l^{-}(w,\phi s_{\alpha})}(\tau^{-}(w, \phi s_{\alpha}),s_{\alpha}
\cdot\mu)$ is surjective. 

On the otherhand, by Lemma 5.4(4), we have $l^{-}(w,\phi)+1=
l^{-}(w,\phi s_{\alpha})$ and $\tau^{-}(w,\phi) s_{\alpha}=
\tau^{-}(w,\phi s_{\alpha})$, and $\tau^{-}(w,\phi s_{\alpha})(\alpha)$ is a 
negative root. Hence,, we have $H^{l^{-}(w,\phi)}(w,\mu)= H^{l^{-}(w,\phi s_{\alpha})}(w,s_{\alpha}\cdot\mu)$, and $H^{l^{-}(w,\phi)}(\tau^{-}(w,\phi s_{\alpha}),\mu)= H^{l^{-}(w,\phi s_{\alpha})}(w,s_{\alpha}\cdot\mu)$, since both  
$w(\alpha)$ and $\tau^{-}(w,\phi s_{\alpha})(\alpha)$ are both negative roots, 
and $\langle\mu, \alpha\rangle$ is positive. Thus, the restriction map $H^{l^{-}(w,\phi)}(w,\mu)\longrightarrow H^{l^{-}(w,\phi)}(\tau^{-}(w, \phi s_{\alpha}),\mu)$ is surjective by the Observation $(1)$. 
Hence, the assertion follows from the surjection 
$H^{l^{-}(w,\phi)}(\tau^{-}(w,\phi s_{\alpha}), \mu)\longrightarrow 
H^{l^{-}(w,\phi)}(\tau^{-}(w,\phi), \mu)\longrightarrow H^{l^{-}(w,\phi)+1}
(\tau^{-}(w,\phi), K)=(0)$, where $K$ is the kernel of the map 
$V_{\mu, \alpha}\longrightarrow \mu$.  

{\bf Case(2):}  $w(\alpha)$ is a positive root.

Observation(2): In this case, $w s_{\alpha}(\alpha)$ is a negative root
and by Lemma 5.4(4), we have $\tau^{-}(w,\phi)=\tau^{-}( w s_{\alpha}, \phi)=\tau^{-}(w s_{\alpha}, \phi s_{\alpha})s_{\alpha}$.
Using the following short exact sequence of $B$- modules:
$(0)\longrightarrow K \longrightarrow V_{\mu, \alpha}\longrightarrow \mu \longrightarrow (0)$, we get the following commutative diagram of $B$- modules:
$$\begin{array}{ccccccccc}
H^{l^{-}(w,\phi)}(w s_{\alpha}, \mu) & \rightarrow & H^{l^{-}(w,\phi)}(w,\mu)& \rightarrow & H^{l^{-}(w,\phi)+1}(w,K)\\  
\downarrow & & \downarrow & & \downarrow \\
H^{l^{-}(w,\phi)}(\tau^{-}(w,\phi)s_{\alpha}, \mu) & \rightarrow & H^{l^{-}(w,\phi)}(\tau^{-}(w,\phi), \mu) & \rightarrow & 
H^{l^{-}(w,\phi)+1}(\tau^{-}(w,\phi), K)=(0)\\  
\end{array}.$$
Since (by Observation(2)), both $w s_{\alpha}(\alpha)$ and 
$\tau^{-}(w s_{\alpha}, \phi s_{\alpha})(\alpha)$ are negative roots,
we have $H^{l^{-}(w,\phi)}(w s_{\alpha}, \mu)= H^{l^{-}(w s_{\alpha},\phi 
s_{\alpha})}(w s_{\alpha},
s_{\alpha}\cdot\mu)$ and 
$H^{l^{-}(w,\phi)}(\tau^{-}(w s_{\alpha},\phi s_{\alpha}), \mu) =H^{l^{-}(w 
s_{\alpha},\phi s_{\alpha})}(\tau^{-}(w s_{\alpha},\phi s_{\alpha}), 
s_{\alpha}\cdot\mu)$.
Since $l(\phi s_{\alpha})=l(\phi)-1$, the left vertical map is surjective,
by induction. The horizontal map $H^{l^{-}(w,\phi)}(\tau^{-}(w,\phi)s_{\alpha}, \mu)\longrightarrow H^{l^{-}(w,\phi)}(\tau^{-}(w,\phi), \mu)$ since cohomology 
vanishes beyond dimension.

{\bf Step(2):} Let $\lambda$ satifies the hypothesis of $(4)$. By taking $\mu:=-(\lambda +\rho)$, we have $\phi\cdot\mu=-\phi(\lambda)-\rho$ and $\langle \phi\cdot\mu,\gamma\rangle \leq 0$ for all simple roots $\gamma$. Hence, by Step(1), $H^{l^{-}(w,\phi)}(w,\mu)$ is nonzero. Now, by Serre-duality, $H^{l(w)-l^{-}(w,\phi)}(w,\lambda-D)=H^{l^{-}(w,\phi)}(w,-\lambda+D-(\rho+D))^{*}=H^{l^{-}(w,\phi)}(w,\mu))^{*}$, and the later is nonzero.
\end{proof}

\bcor Let $w\in W$ be such that $s_{\alpha}\leq w$ for all simple roots $\alpha$. If $\lambda$ is a weight such that $w\cdot \lambda$ is dominant and with the property that $\langle w(\lambda), \gamma \rangle \geq -1$ for all simple roots $\gamma$, then, either $G$ is finite dimensional and $X(w)=G/B$ or there are two integers $i ~ < ~ j$ such that set $H^{i}(w, \lambda - D)\neq 0$ and $H^{j}(w,\lambda)\neq 0$.\ecor
Proof of this corollary follows from Theorem 5.8(1) and 5.8(4).

\brem Even though our proof of Theorem 5.8(2) holds only in the finite dimensional case, we believe that the statement holds even in the infinite dimensional case. \erem  

\brem Although Theorem 5.8(4) is stated with a twist by the boundary divisor, we believe that the $H^{l(w)-l^{-}(w,\phi)}(w,\lambda)$ itself is nonzero for any generic $\lambda$ in the $\phi$-chamber. \erem   

\brem We believe that for a generic weight $\lambda$ in the $\phi$-chamber the $B$- module $H^{l^{+}(w,\phi)}(w,\lambda)$ is a cyclic $B$- module generated by a weight vector of weight $\tau^{+}(w,\phi)\cdot\lambda$. \erem

\section{A Surjectivity Theorem} In this section, we prove that if $G$ is finite dimensional and $\lambda$ is a generic weight in the $\phi$-chamber and $\phi\leq w$, then, the restriction $H^{l(\phi)}(G/B,\lambda)\longrightarrow H^{l(\phi)}(w,\lambda)$ is surjective.

Let $\phi \in W$. Let $M:=max\{\langle \beta , \gamma \rangle : \beta\in \phi(S), \gamma\in S$.
Let $\lambda$ be a weight such that $\langle \phi\cdot\lambda, \gamma \rangle > l(\phi)M$ for all $\gamma\in S$. Let $w\in W$ be such that $\phi \leq w$. Then, we have: \bth The restriction map $H^{l(\phi)}(G/B, \lambda)\longrightarrow 
H^{l(\phi)}(w, \lambda)$ is surjective. \eeth
\begin{proof} Proof is by induction on $l(\phi)$.
If $l(\phi)=0$, then, $\lambda$ is dominant, and so there is nothing to prove. 
So, let $\phi$ be such that $l(\phi)$ is positive.
Let $w\in W$ and let $\phi\leq w$. 
Let $\lambda$ be as in the hypothesis above with respect to $\phi$. 
Choose an $\alpha\in S$ such that $\phi(\alpha)<0$. For this choice of 
$\alpha$, we have two possibilities:

{\bf Case(1):} $w(\alpha)< 0$. We have $H^{l(\phi)}(w, \lambda)=H^{l(\phi)-1}(w, s_{\alpha}.\lambda)$. Since $\langle \phi s_{\alpha}\cdot s_{\alpha}\cdot\lambda, \gamma \rangle=\langle \phi\cdot \lambda, \gamma \rangle$ and this number is bigger than $l(\phi)M = (1+l(\phi s_{\alpha}))M$, by induction, the restriction map $H^{l(\phi s_{\alpha})}(G/B, s_{\alpha}\cdot \lambda)\longrightarrow 
H^{l(\phi s_{\alpha})}(w, s_{\alpha}\cdot \lambda)$ is surjective, and this map is the same as $H^{l(\phi)}(G/B, \lambda)\longrightarrow H^{l(\phi)}(w,\lambda)$, and we are done.

{\bf Case(2):}  $w(\alpha) > 0$. 
Using the following short exact sequence of $B$- modules
$(0)\longrightarrow \lambda\longrightarrow V_{s_{\alpha}(\lambda), \alpha}
\longrightarrow Q\longrightarrow (0)$, we have exact sequence of $B$-modules:
$H^{l(\phi s_{\alpha})}(w, Q)\longrightarrow H^{l(\phi) }(w, \lambda)
\longrightarrow H^{l(\phi)}(ws_{\alpha}, s_{\alpha}(\lambda))$.
From the hypothesis of $\lambda$, we have $\langle \phi s_{\alpha}\cdot
s_{\alpha}(\lambda), \gamma \rangle ~ > ~ l(\phi s_{\alpha})M$ for all 
simple roots $\gamma$. Therefore, by Proposition 5.1, 
$H^{l(\phi)}(w s_{\alpha}, s_{\alpha}(\lambda))=(0)$, since 
$l(\phi) ~ > ~ l(\phi s_{\alpha})$. 
Thus, we have 

{\bf Observation (1):} The natural map 
$H^{l(\phi s_{\alpha})}(w, Q)\longrightarrow 
H^{l(\phi)}(w,\lambda )$ is surjective.  

Now, by using the $P_{\alpha}/B$- fibration $G/B\longrightarrow G/P_{\alpha}$,
we have $H^{l(\phi)}(G/B, V_{s_{\alpha}(\lambda), \alpha})=H^{l(\phi)}(G/B, 
s_{\alpha}(\lambda))$. Thus, we have $H^{l(\phi)}(G/B, V_{s_{\alpha}(\lambda), \alpha})=H^{l(\phi)}(G/B, s_{\alpha}(\lambda))= (0)$. Hence, using the exact 
sequence of $B$- modules $H^{l(\phi s_{\alpha})}(G/B, Q)\longrightarrow H^{l(\phi) }(G/B, \lambda)\longrightarrow H^{l(\phi)}(G/B, V_{s_{\alpha}(\lambda), \alpha})=H^{l(\phi)}(G/B, s_{\alpha}(\lambda))=(0)$ since $l(\phi) > l(\phi s_{\alpha})$. Thus, we have

{\bf Observation(2):} The natural map 
$H^{l(\phi s_{\alpha})}(G/B, Q)\longrightarrow H^{l(\phi)}(G/B, \lambda)$
is surjective. 
Now, using the short exact sequence of $B$- modules, 
$(0)\longrightarrow V_{s_{\alpha}\cdot \lambda, \alpha}\longrightarrow Q 
\longrightarrow s_{\alpha}(\lambda)\longrightarrow (0)$, 
we have the two following exact sequences of $B$- modules:
$$H^{l(\phi s_{\alpha})}(w s_{\alpha}, s_{\alpha}\cdot \lambda)\longrightarrow 
H^{l(\phi s_{\alpha})}(w, Q)\longrightarrow H^{l(\phi s_{\alpha})}(w, 
s_{\alpha}(\lambda))\longrightarrow H^{l(\phi)}(w s_{\alpha}, s_{\alpha}
\cdot\lambda)=(0),$$ and $H^{l(\phi s_{\alpha})}(G/B, s_{\alpha}\cdot\lambda)=H^{l(\phi s_{\alpha})}(G/B, V_{s_{\alpha}\cdot\lambda, \alpha})\longrightarrow 
H^{l(\phi s_{\alpha})}(G/B, Q)\longrightarrow H^{l(\phi s_{\alpha})}(G/B, 
s_{\alpha}(\lambda))=H^{l(\phi s_{\alpha})}(G/B, V_{s_{\alpha}(\lambda), 
\alpha})\longrightarrow H^{1+l(\phi s_{\alpha})}(G/B, s_{\alpha}\cdot\lambda)
=(0)$. 

The isomorphisms $H^{i}(G/B, s_{\alpha}\cdot\lambda)=H^{i}(G/B, V_{s_{\alpha}\cdot\lambda, \alpha})$, for $i=l(\phi s_{\alpha}), l(\phi)$ and $H^{l(\phi s_{\alpha})}(G/B, s_{\alpha}(\lambda))
=H^{l(\phi s_{\alpha})}(G/B, V_{s_{\alpha}(\lambda), \alpha})$ follows from a similar argument using the $P_{\alpha}/B$-fibration $G/B\longrightarrow G/P_{\alpha}$. The ismorphisms $H^{1+l(\phi s_{\alpha})}(G/B, s_{\alpha}\cdot\lambda)=(0)$ follows from Bott's Theorem. 

Let $V_{1}$ be the image of $H^{l(\phi s_{\alpha})}(w s_{\alpha},  s_{\alpha}
\cdot \lambda)\longrightarrow H^{l(\phi s_{\alpha})}(w, Q)$, 
$V=H^{l(\phi s_{\alpha})}(w, Q)$, and let $V_{2}=H^{l(\phi s_{\alpha})}(w, 
s_{\alpha}(\lambda))$.  Also, let $V_{1}^{\prime}$ be the image of 
$H^{l(\phi s_{\alpha})}(G/B, s_{\alpha}\cdot\lambda)=H^{l(\phi s_{\alpha})}
(G/B, V_{s_{\alpha}\cdot\lambda, \alpha})\longrightarrow 
H^{l(\phi s_{\alpha})}(G/B, Q)$, and let $V^{\prime}=H^{l(\phi s_{\alpha})}
(G/B, Q)$, and let $V_{2}^{\prime}= H^{l(\phi s_{\alpha})}(G/B, 
s_{\alpha}(\lambda))$. We therefore have the following commutative diagram of two short exact sequences of $B$- modules: $$\begin{array}{ccccccccc}
(0) & \longrightarrow & V_{1}^{\prime} & \longrightarrow &
V^{\prime} & \longrightarrow & V_{2}^{\prime} &
\longrightarrow & (0) \\
\downarrow & & \downarrow & & \downarrow & & \downarrow & & \downarrow \\
(0) & \longrightarrow & V_{1} & \longrightarrow & V & \longrightarrow & V_{2} & \longrightarrow & (0). \\
\end{array}$$

{\bf Step(1):} We now claim that the first vertical map 
$V_{1}^{\prime}\longrightarrow V_{1}$ and the third vertical map
$V_{2}^{\prime}\longrightarrow V_{2}$ are surjective.  
Now, $\langle \phi s_{\alpha}\cdot s_{\alpha}(\lambda), \gamma \rangle 
= \langle \phi\cdot\lambda -\phi(\alpha), \gamma \rangle > (l(\phi)-1) M 
=l(\phi s_{\alpha})M$, and $\langle\phi s_{\alpha}\cdot s_{\alpha}\cdot\lambda,
\gamma \rangle = \langle \phi\cdot\lambda, \gamma \rangle > l(\phi)M$ for all 
simple roots  $\gamma$. Therefore, by induction, the maps 
$H^{l(\phi s_{\alpha})}(G/B, s_{\alpha}\cdot \lambda)\longrightarrow 
H^{l(\phi s_{\alpha})}(w s_{\alpha}, s_{\alpha}\cdot\lambda)$ 
and $H^{l(\phi s_{\alpha})}(G/B, s_{\alpha}(\lambda))\longrightarrow 
H^{l(\phi s_{\alpha})}(w s_{\alpha}, s_{\alpha}(\lambda))$ are surjective. 
 Hence, the proof of Step(1) follows from the definitions of $V_{1}$, 
$V_{1}^{\prime}$, $V_{2}$ and $V_{2}^{\prime}$.  

{\bf Step(2):} The middle vertical map $V^{\prime}\longrightarrow V$ is surjective. The proof of Step(2) follows from Step(1) and snake Lemma. 

{\bf Step(3):} We now complete the proof of the Theorem.
Consider the commutative diagram of $B$- modules :
$$\begin{array}{ccccccccc}
H^{l(\phi s_{\alpha})}(G/B, Q) & \longrightarrow &
H^{l(\phi)}(G/B, \lambda) & \longrightarrow & (0)\\
\downarrow & & \downarrow & & \downarrow \\
 H^{l(\phi s_{\alpha})}(w, Q) & \longrightarrow & H^{l(\phi)}(w, \lambda) 
& \longrightarrow & (0). \\ \end{array}$$
Here $Q$ is the $B$- module defined above. The surjectivity of the second horizontal map follows from Observation (1). The surjectivity of the left vertical map follows from Step(2). So, the second vertical map is also surjective, completing the proof
\end{proof}

Chennai Mathematical Institute, 92, G.N.Chetty Road, T.Nagar, Chennai
600 017, India
(email: kannan@cmi.ac.in).

\end{document}